\documentclass[default,iicol]{sn-jnl}% Default with double column layout

\usepackage{amsmath, amsthm, amssymb, graphicx}

\jyear{2022}%

\theoremstyle{thmstyleone}%
\newtheorem{theorem}{Theorem}%  meant for continuous numbers

\newtheorem{alg}[theorem]{Algorithm}
\newtheorem{lemma}[theorem]{Lemma}
\newtheorem{example}[theorem]{Example}%
\newtheorem{remark}[theorem]{Remark}%
\newtheorem{definition}[theorem]{Definition}%

\begin{document}

\title[~]{\huge Application of integral invariants to apictorial jigsaw puzzle assembly}


\author[1]{\fnm{Peter} \sur{Illig}}%\email{pkillig@gmail.com}

\author[2]{\fnm{Robert} \sur{Thompson}} 
\equalcont{Corresponding author:  \url{rthompson@carleton.edu}}

\author[3]{\fnm{Qimeng} \sur{Yu}}%\email{QimengYu2023@u.northwestern.edu}

\affil[1]{\orgdiv{Epic Systems Corporation}}

\affil[2]{\orgdiv{Department of Mathematics and Statistics}, \orgname{Carleton College}}

\affil[3]{\orgdiv{McCormick School of Engineering}, \orgname{Northwestern University}}

%%==================================%%
%% Abstract
%%==================================%%

\abstract{We present a method for the automatic assembly of apictorial jigsaw puzzles.  This method relies on integral area invariants for shape matching and an optimization process to aggregate shape matches into a final puzzle assembly.  Assumptions about individual piece shape or arrangement are not necessary.   We illustrate our method by solving example puzzles of various shapes and sizes.}




\maketitle

\section{Introduction}

We present a method for the automatic assembly of jigsaw puzzles.  Our approach is apictorial, using only shape information provided by the boundaries of the pieces.  The method has three basic components:
\begin{enumerate}
\item Compute an {\it integral invariant} for each puzzle piece, encoding its shape independent of position and orientation.
\item Compare integral invariants to determine matches among pairs of pieces and measure the quality of these matches.
\item Assemble the puzzle by aggregating these pairwise matches as consistently as possible.
\end{enumerate}

The development of computational approaches to the solution of jigsaw puzzles using only shape information began as early as 1964, \cite{freeman1964apictorial}.  Much of the ensuing work, e.g. \cite{wolfson1988solving, kosiba1994automatic, goldberg2002global}, has focused on traditional rectangular jigsaw puzzles, and leverages assumptions about piece shape and puzzle arrangement in the solving process.  Key assumptions among these are that puzzle pieces are four-sided with ``indents'' and ``outdents'' and that there are corner and edge pieces which may be identified and assembled separately from interior pieces.  A notable exception is \cite{hoff2014automatic}, where an extended method of differential invariant signatures, \cite{hoff2013extensions}, and an intensive piece locking method is utilized that can be effective in assembling both standard (rectangular) and nonstandard puzzles, \cite{rainforest, baffler}, without these usual assumptions.  The present work proceeds in the spirit of \cite{hoff2014automatic, grim2016automatic}, eschewing structural information about the pieces and their arrangement.

As observed in \cite{goldberg2002global}, apictorial jigsaw puzzle assembly has two main difficulties:  the \textit{geometric difficulty} of reliably determining when pieces fit together, and the \textit{combinatorial difficulty} of parsing the very large number of ways that the collection of pieces can be assembled.  To address the geometric difficulty we apply a matching process which attempts to find the longest fits within a threshold of shape similarity, using integral area invariants to measure this similarity.  We find that this process can correctly identify entire matching sides of rectangular pieces and provide visually close fits in a wide variety of examples, eliminating the necessity of a piece-locking process like that of \cite{hoff2014automatic}.  To address the combinatorial difficulty we characterize a puzzle assembly as the solution to an optimization process, measuring the cost of including a particular piece fit in this optimal assembly via a combination of a  local  measurement (based on quantities computed from the fit itself) and a global  measurement (based on consistency of collections of piece fits).

\begin{figure}[htbp] %  figure placement: here, top, bottom, or page
   \centering
   \includegraphics[width= \linewidth]{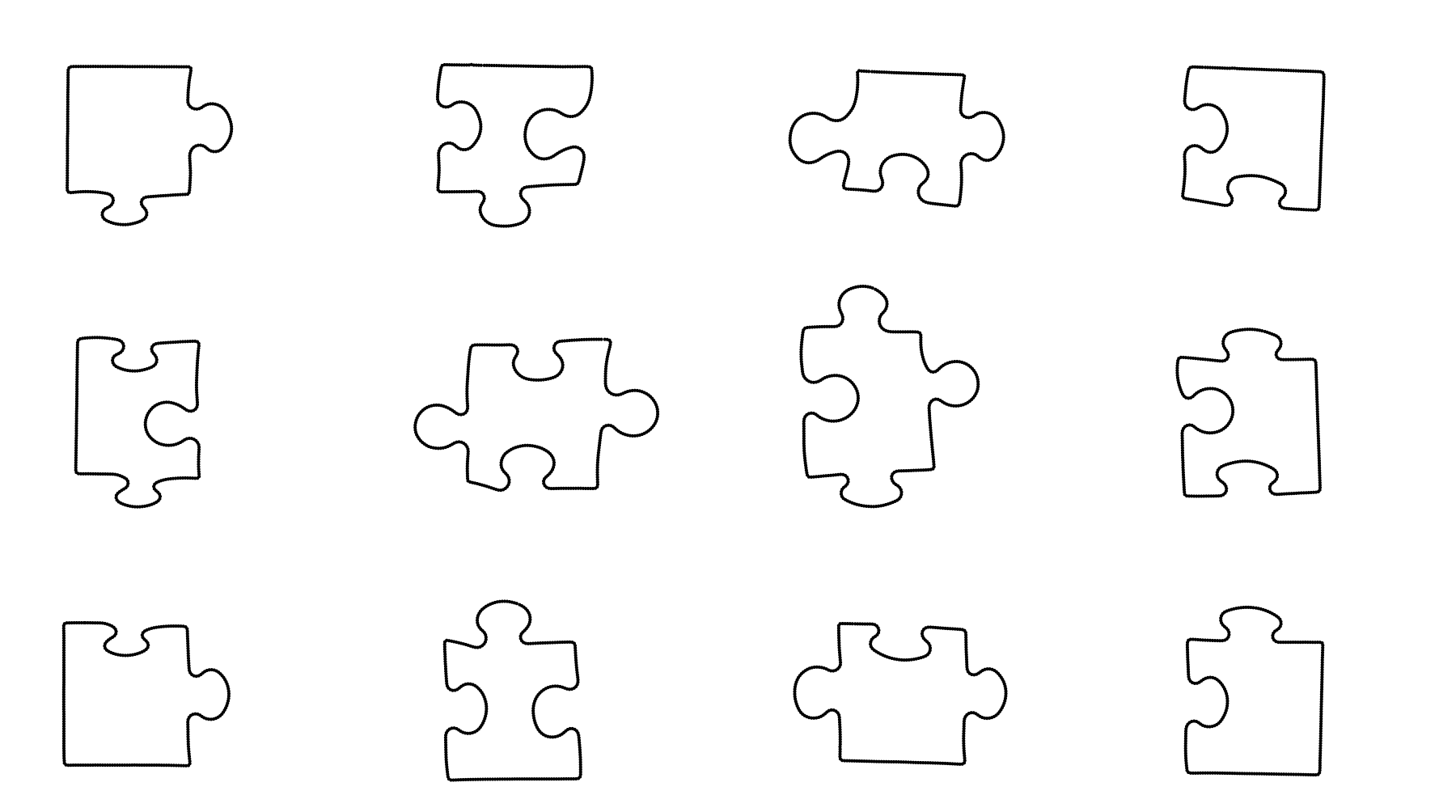} 
   \caption{The 12 piece puzzle used for examples throughout the paper.}
   \label{fig:examplepuzzle}
\end{figure}

The layout of the paper is as follows.  Section \ref{sec:data} describes the acquisition and preprocessing of the jigsaw puzzle data and outlines in Algorithm \ref{alg:resampling} a method for respacing the discrete curve information that produces more reliable comparisons of curve shapes.  Section \ref{sec:invariants} recalls the definition of the integral area invariant, and provides in Algorithm \ref{alg:areainvariant} a way to compute this invariant exactly for piecewise linear curves.  Section \ref{sec:comparison}  addresses piece comparison and the alignment of integral invariants used for finding the ``best'' fit between a pair of puzzle pieces.  In Section \ref{sec:assembly} the piece fits are aggregated into a puzzle assembly, and the criterion of cycle consistency is introduced to ensure compatibility of collections of piece fits.  Finally, we illustrate our algorithmic approach through various examples in Section \ref{sec:results}.  We motivate and demonstrate our methods throughout the paper on a simple 12 piece example puzzle, \cite{dino}, shown in Figure \ref{fig:examplepuzzle}.

\section{Puzzle data}
\label{sec:data}
\begin{figure}[htbp] %  figure placement: here, top, bottom, or page
   \centering
   \includegraphics[width=.6 \linewidth]{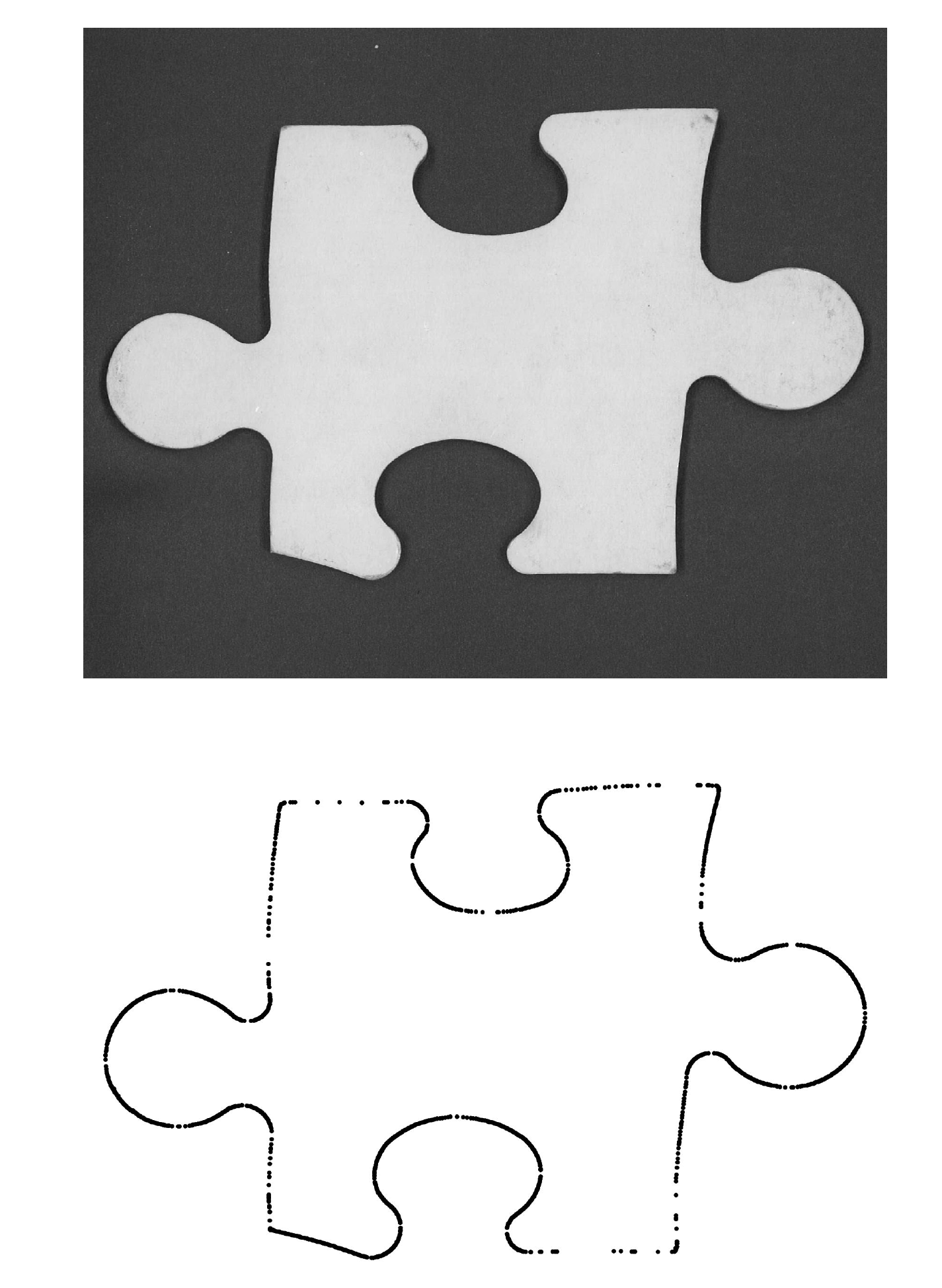} 
   \caption{A puzzle piece image and its (unprocessed) boundary curve.}
   \label{fig:segmentation}
\end{figure}

The input to our algorithm is a collection of ordered lists of points.  Each element of the collection represents a puzzle piece, and each ordered list is a sampling of points around the boundary of the puzzle piece.  To obtain this sampling of the boundary, each puzzle piece is photographed via a photocopier, then processed using image segmentation.  In our computations, segmentation was done in \texttt{Mathematica} by using \texttt{MorphologicalBinarize} to create a binary image, and \texttt{ComponentMeasurements} to extract the boundary after binarization.  The same task could be accomplished in \texttt{Matlab} using the command \texttt{bwboundaries}, or in other software using edge detection or active contour methods, \cite{kass1988snakes}.  An example puzzle piece image and its segmented boundary are shown in Figure \ref{fig:segmentation}.

The raw output of boundary points from segmentation can be unevenly spaced.  Since our assembly method relies on comparing shape signatures as a function of arclength, this raw output must be resampled so that consecutive boundary points are separated by the same fixed arclength.  This is accomplished via repeated linear interpolation and resampling according to a fixed arclength measurement.  The idea of this method is suggested in \cite{hoff2014automatic} and studied more carefully in \cite{manivel2021iterative}, where it is shown that this repeated interpolation will indeed converge to an evenly spaced discrete curve.   This method is a fundamental step in preprocessing the data for shape comparison, and we outline it in more detail in Algorithm \ref{alg:resampling}.

%%%%%%%%%%%%%%%%%%%%%

\begin{alg} Resampling a closed discrete curve by a fixed arclength.
 
  \label{alg:resampling}
  \vskip 5pt
  \noindent
 \textbf{Input:}  An ordered collection of points $p_0, \ldots, p_n$ in $\mathbb{R}^2$, with $p_0 = p_n$, representing a sampling of a closed curve.  Adjacent points should be distinct:  $p_k \neq p_{k+1}$.
  \vskip 5pt
  \noindent
  \textbf{Output:}  An ordered collection of points $q_0, \ldots, q_m$ in $\mathbb{R}^2$ representing a new sampling of the closed curve satisfying $\vert\vert q_{k+1} - q_{k}\vert\vert = \delta$, $k = 0, \ldots, m-1$ for some chosen distance $\delta$.
   \vskip 5pt
  \begin{enumerate}
    \item  Let $d_0 = 0$ and recursively compute $d_k = \vert\vert p_k - p_{k-1}\vert\vert + d_{k-1}$ for $k=1, \ldots n$.  $d_k$ is the piecewise linear arclength distance from $p_0$ to $p_k$.

    \item
    Compute the piecewise linear interpolating function $g:[0,d_n] \rightarrow [0,n]$ for the points $(d_k,k)$, $0 \leq k \leq n$.  This function inverts the arclength measurements, so that $g(d_k) = k$.  Here we require $p_k \neq p_{k+1}$ in order for this inverse to be well defined.
    
   \item
  Compute the piecewise linear interpolating function $h:[0,n] \rightarrow \mathbb{R}^2$ for the discrete curve points $p_0, \ldots, p_n$.
  
  \item  Choose a fixed arclength $\delta$ and compute a new collection of points $q_0, \ldots, q_m$ via $q_k = h(g(k \delta))$ for  $k = 0, 1, \ldots, m= \lfloor d_n/\delta \rfloor$.  This is a new sampling of the discrete curve, where points are separated by a distance of $\delta$.
  
  \item  Optionally, set $q_{m+1} = q_0$ and repeat steps 1-4 using the new collection $q_0, \ldots, q_{m+1}$ to further smooth the discrete curve.
  
  \end{enumerate}
\end{alg}
%%%%%%%%%%%%%%%%%%%%%

\begin{figure}[h] 
   \centering
   \includegraphics[width=.95 \linewidth]{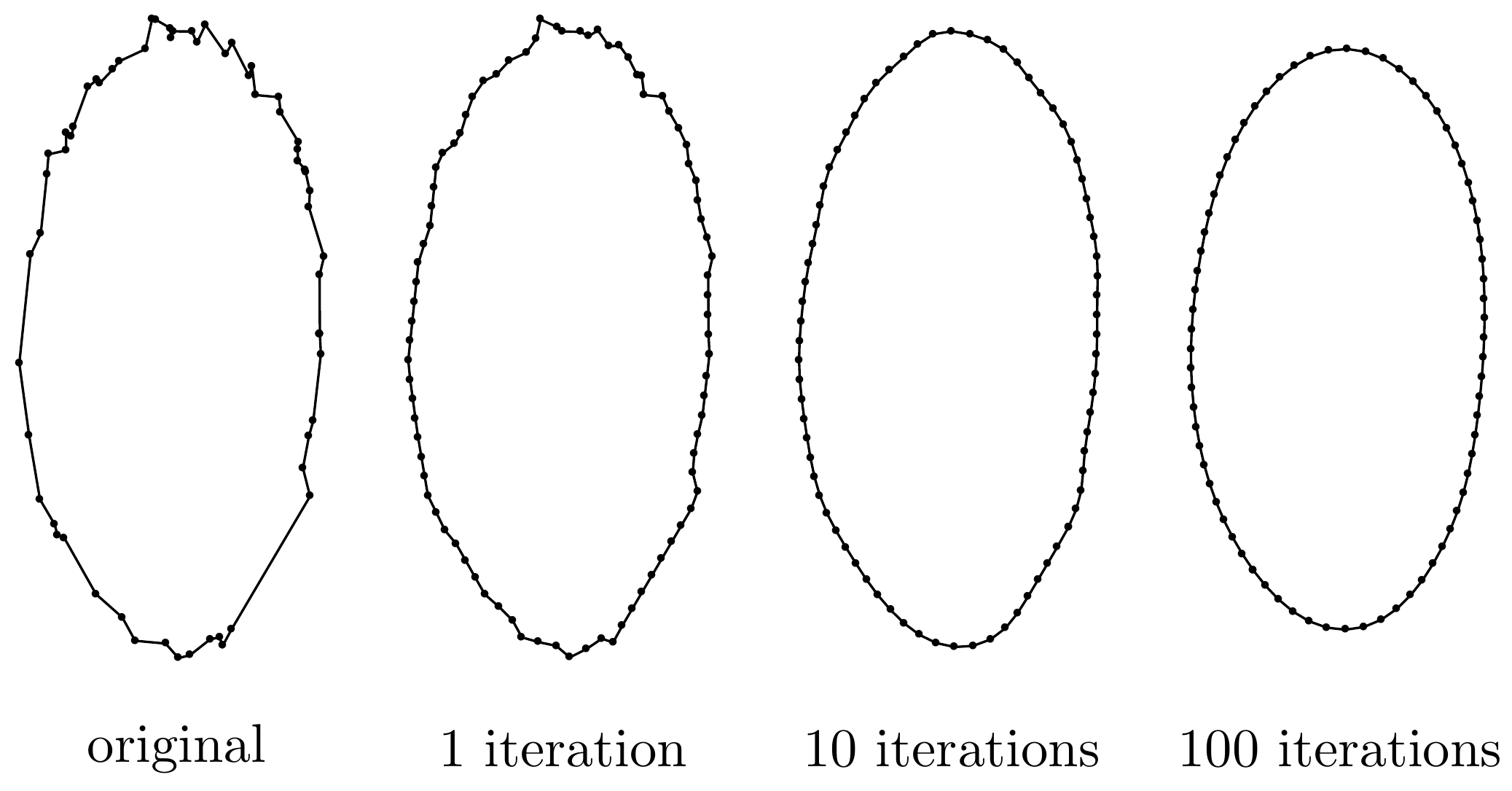} 
   \caption{A noisy discrete curve and results after iterations of Algorithm \ref{alg:resampling}.}
   \label{fig:noisyellipse}
\end{figure}

\begin{example}
We illustrate empirically the effect of applying Algorithm \ref{alg:resampling} to a curve with widely varying distances between points and added noise.  This illustration is shown in Figure   \ref{fig:noisyellipse}. We choose 60 random points on an ellipse, and add a small amount of noise in the radial direction.  After 1 iteration, the points are nearly uniformly spaced by arclength.  Further iteration continues to smooth and shorten the curve.  This shortening is an artifact of Step 4 in Algorithm \ref{alg:resampling}, where a small ``leftover'' part of the curve is discarded.
\end{example}

\begin{remark}
Because our goal is direct comparison of puzzle piece boundaries, we sample all boundaries using the same arclength, typically leaving a single anomalous distance $\vert\vert q_{m} - q_0\vert\vert \neq \delta$ after applying Algorithm \ref{alg:resampling}.  We did not encounter any issues arising from this anomolous distance.
\end{remark}

\begin{figure}[h] 
   \centering
   \includegraphics[width=.6 \linewidth]{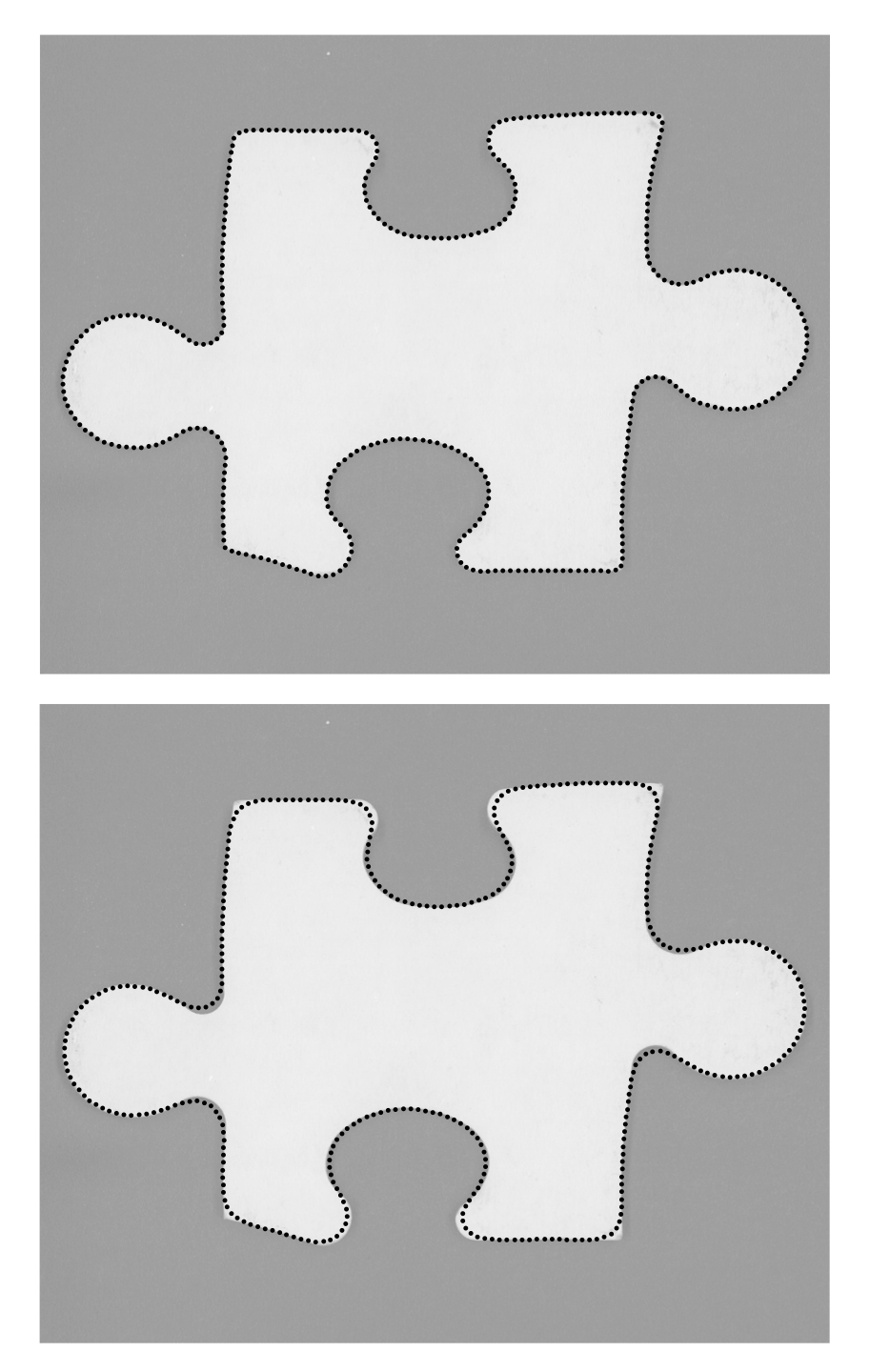} 
   \caption{A puzzle boundary with $\delta = 20$ after $5$ and $30$ smoothing iterations.}
   \label{fig:boundarysmoothing}
\end{figure}

After obtaining the unprocessed puzzle piece boundaries via segmentation, Algorithm \ref{alg:resampling} is applied to each boundary for a predetermined number of iterations and arclength distance $\delta$.  The distance $\delta$ is chosen to balance the precision of pairwise comparison of pieces with the computational time needed for comparison.  Using a pixel's length or width for the unit distance, values of $\delta$ in the range of $5$ to $20$ performed well for all examples, based on puzzle piece images with a resolution of 300 pixels per inch.  The number of iterations of Algorithm \ref{alg:resampling}  to achieve a visually acceptable balance of smoothness and accuracy varied from $5$ to $30$, depending on the distance $\delta$.   Figure \ref{fig:boundarysmoothing} illustrates this visual selection process at a resolution of $\delta = 20$;  $5$ iterations results in an accurate representation of the piece boundary, while $30$ iterations shows excessive smoothing.

\section{Integral area invariants} 
\label{sec:invariants}

To determine if two digital jigsaw puzzle pieces fit together, we  compare the shapes of their boundaries, searching for portions of each boundary that are congruent under some rotation and translation in the plane (the action of some element of the special Euclidean group $SE(2)$).  Invariants facilitate this shape comparison by removing the freedom of  rotation and translation; puzzle matches can be found via direct comparison of the invariants rather than the pieces themselves.  We focus on a simple integral invariant, the integral area invariant,  first introduced in \cite{manay2004integral} and studied for its shape identification properties, \cite{manay2006integral}.

Let $p:S^1 \rightarrow \mathbb{R}^2$ be a closed simple planar curve,  $R$ the region enclosed by $p$, and  $\partial R$ the boundary of $R$ (which is also the image of $p$).

\begin{figure*}[htbp] %  figure placement: here, top, bottom, or page
   \centering
   \includegraphics[width=.7\textwidth]{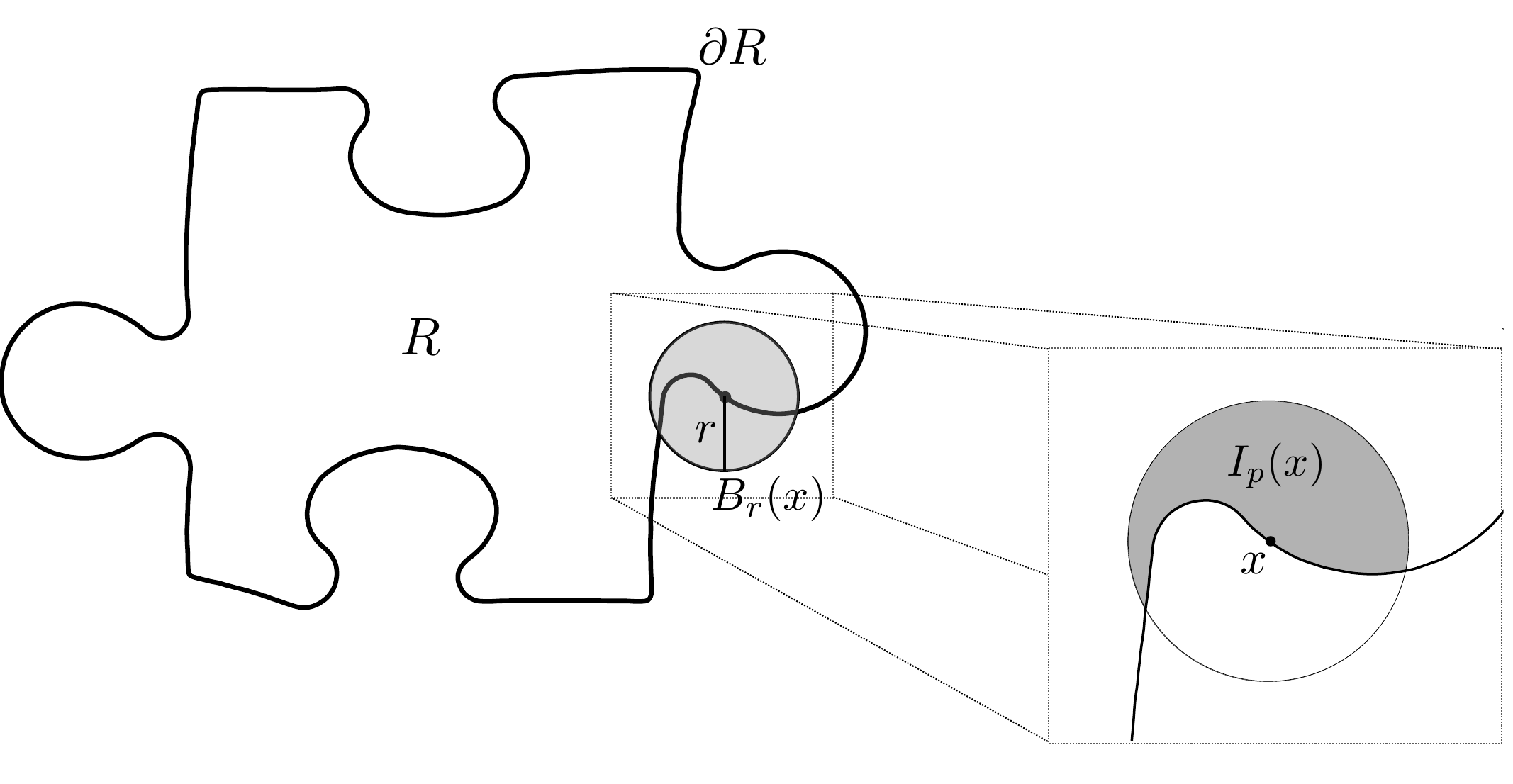} 
   \caption{Visualization of the integral area invariant at a point $x$.}
   \label{fig:integralinvariant}
\end{figure*}

\begin{definition} \label{def:integralarea}
Let $r>0$ and let $B_r(x)$ be the disk of radius $r$ centered at $x$.  The \textit{integral area invariant} (of radius $r$) for $p$ is given by
\[
I_p(x) =   \int_{B_r(x) \cap R} dA,
\]
 the area of the intersection of $B_r(x)$ and $R$.  This definition is illustrated in Figure \ref{fig:integralinvariant}.
\end{definition}

The integral area invariant is invariant under the action of $SE(2)$.  This means that
\[
I_p(x) = I_{g p}( g x) \quad   \text{for all } g \in SE(2),
\]
where $g p$ is the curve $p$ transformed by $g$.  By virtue of this $SE(2)$ invariance, the integral area invariant provides shape information -- information about the image of $p$ which is independent of  placement in the plane.  Congruent curves must have the same integral area invariants, but the extent to which an integral area invariant uniquely determines the  curve up to congruence is the subject of ongoing research, \cite{fidler2007inverse, fidler2008identifiability, calder2012circular}.  It is enough for us to know that, practically speaking, integral invariants will help us identify when two curves have very similar shape.

The integral area invariant makes sense for any $x$ in $\mathbb{R}^2$, but in practice, we will restrict the domain to the image of $p$.  With a parameterization of $p$ in hand, $I_p$ can be interpreted as a function of the parameter; for the curve $p(s)$ in $\mathbb{R}^2$ we obtain the real valued function $I_p(s)$.

To use integral invariants for puzzle piece comparison,  we  adapt the above discussion to the discrete setting.  In the following, let $p_0, \ldots, p_n$ be a collection of points in $\mathbb{R}^2$ representing a closed curve without self-intersections, and  $p_0 = p_n$.  As before, these points represent a sampling of the outline of a puzzle piece.  The \textit{discrete integral area invariant} is then defined at the points $p_0, \ldots, p_n$ just as in Definition \ref{def:integralarea} by taking $p$ to be the piecewise linear interpolation of $p_0, \ldots, p_n$.

We now outline an exact method for computing this discrete version of the integral area invariant.  See also \cite{o2020computation}.  We will need to find the area of a large polygon, for which we will use the following classical formula, sometimes called the shoelace formula.

\begin{lemma} \label{shoelace}
Suppose that  $q_0, \ldots, q_n$ are the sequential vertices of an $n$ sided polygon, with $q_0 = q_n$.  The area of this polygon is  given by the formula
\[
\frac{1}{2}  \left \vert \sum_{k = 1}^n \det \bigg ( q_{k-1}  \,\,\,\, q_k \bigg )   \right \vert.
\]
 \end{lemma}

\begin{figure}[htbp] %  figure placement: here, top, bottom, or page
   \centering
  \includegraphics[width=.8\linewidth]{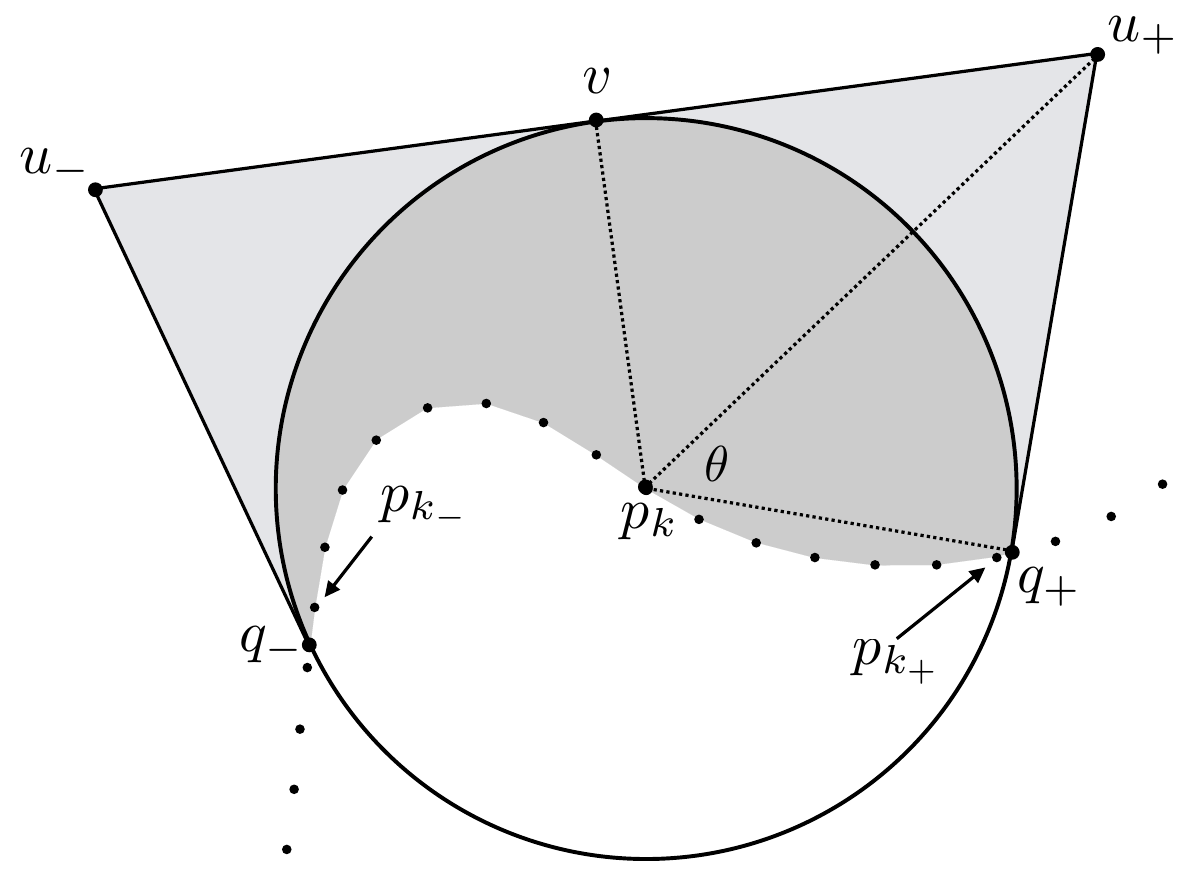} 
   \caption{Constructions for computing the discrete integral area invariant.}
   \label{fig:discreteintegralarea}
\end{figure}

In the following discussion we take indices modulo $n$ and assume that the points outline the curve in counter-clockwise orientation.  Refer to Figure \ref{fig:discreteintegralarea} for an illustrations of the calculations.
%%%%%%%%%%%%%%%%%%%%%
\begin{alg} Calculating the discrete integral area invariant.
 \label{alg:areainvariant}
  \vskip 5pt
  \noindent
 \textbf{Input:}  A fixed radius $r$, and an ordered collection of points $p_1, \ldots, p_n$ in $\mathbb{R}^2$ representing a closed curve without self-intersections.  
  \vskip 5pt
  \noindent
  \textbf{Output:}  The integral area invariant value $I_p(p_k)$ where $p$ is the piecewise linear interpolation of the points $p_0, \ldots, p_n$.  
   \vskip 5pt
  \begin{enumerate}
    \item  Let $B_r(p_k)$ be the disk of radius $r$ centered at $p_k$.  Choose $p_{k_+}$ (respectively $p_{k_-}$) to be the final point in the list $p_k, p_{k+1}, \ldots $ (respectively $p_k, p_{k-1}, \ldots$) contained in $B_r(p_k)$.
  
   \item Let $t_\pm$ be the positive solutions to $\vert\vert p_k + t (p_{k_\pm \pm 1} -p_{k_\pm}) \vert\vert^2 = r^2$.   Next, let $q_\pm = p_k + t_\pm (p_{k_\pm \pm 1}-p_{k_\pm})$.  The point $q_+$ (respectively $q_-$) is the intersection of the line between $p_{k_+}$ and $p_{k_++1}$ (respectively $p_{k_-}$ and $p_{k_--1}$) and the boundary of $B_r(p_k)$.
   
   % so
  % \[
  % t_\pm = \frac{-2 p_k \cdot (p_{k_\pm \pm 1}-p_{k_\pm}) + \sqrt{\big( p_k \cdot (p_{k_\pm \pm 1}-p_{k_\pm})\big)^2 - 4 \vert\vert p_{k_\pm \pm 1}-p_{k_\pm}  \vert\vert^2 \big (
  % \vert\vert p_k\vert\vert^2 - r^2\big )}}{2 \vert \vert p_{k_\pm \pm 1}-p_{k_\pm}  \vert \vert^2}.
 %  \]

   \item  Construct points $u_-, u_+, v$ as follows.  Let $\theta = \frac{1}{4} \angle q_+ p_k q_- $, where  $\angle q_+ p_k q_-$ is the angle between $q_+-p_k$ and $q_- - p_k$, measured counter-clockwise.  Define 
   \[
   u_\pm = p_k + \sec( \theta) R_{\pm \theta} (q_\pm- p_k),
   \] 
   where $R_\theta$ is the standard (counter-clockwise) rotation matrix through angle $\theta$.  Then, define $v$ to be the bisector of the segment connecting $u_+$ and $u_-$, or equivalently
  \[
  v = p_k + R_{2 \theta} (q_+ - p_k).
  \]
 
  \item   
  Let $A_1$ be the area of the polygon $u_-, q_-,p_{k_-}, \ldots p_k, \ldots p_{k_+},q_+,u_+$ (which is easily computed via Lemma \ref{shoelace}).  This is the combined area of the light and dark shaded regions in Figure \ref{fig:discreteintegralarea}.

  \item  Let $A_2 = \vert\vert(v-u_+) \times (q_+ - u_+) \vert\vert - r^2 \left ( \frac{\theta}{2} - \sin \frac{\theta}{2}\right)$, the total area of the lighter shaded regions in Figure \ref{fig:discreteintegralarea}.  (The first term is twice the area of the triangle formed by $v, q_+, u_+$, and the second term is twice the area of the circular segment cut by the secant connecting $q_+$ and $v$.  Both terms are doubled by symmetry of the computation.)
  
  \item  Return $I_p(p_k) = A_1 - A_2,$ the discrete integral invariant at $p_k$.

\end{enumerate}
\end{alg}
%%%%%%%%%%%%%%%%%%%%%

\begin{remark}
In Algorithm \ref{alg:areainvariant} what is actually computed is different than $I_p$ as given in Definition \ref{def:integralarea}.  It is assumed that curve does not ``wander back into $B_r(x)$'', so that the intersection of $B_r(x)$ with the exterior of $R$ has only a single component.  For a large enough choice of $r$, this assumption would generally be false.  But, for the purposes of puzzle piece matching,  we are interested in matching portions of curves, not global shape matching, and so prefer $I_p$ from Algorithm \ref{alg:areainvariant} to that of Definition \ref{def:integralarea}.
\end{remark}

\section{Piece comparison} 
\label{sec:comparison}

We now describe the process of comparing two puzzle pieces to find possible fits.  Integral area invariants reduce this problem to one of local sequence alignment; a match between two puzzle pieces is a partial overlap of the integral area invariant sequences, and searching for the ``best'' partial overlap should typically reveal the ``best'' fit possible for a given pair of puzzle pieces.  The meaning of ``best'' here is qualitative and will be illustrated via examples; in particular we will see ways in which good signature overlap may result in poor puzzle piece matches.  

\subsection{Signature alignment and piece fitting}
Before discussing puzzle pieces, we describe for context the process of finding alignments of two periodic arrays to within a fixed threshold.  Let $A = \{a_0, \ldots, a_{m-1}\}$ and $B = \{b_0, \ldots, b_{n-1}\}$ be numerical arrays.  $A$ and $B$ are considered as periodic arrays, so indices are taken modulo their respective array lengths:  if $j, k$ are integers, let $a_j = a_{j\, \text{mod}\, m}$ and $b_k = b_{k\, \text{mod}\, n}$, where representatives are chosen so that $0 \leq j\, \text{mod}\, m \leq m-1$ and $0 \leq k\, \text{mod}\, n \leq n-1$.  This periodicity convention will be used for the remainder of Section \ref{sec:comparison}.  The following definition makes the notion of alignment precise.

\begin{definition}
Let $\epsilon > 0$.  A pair of substrings $\{a_i, a_{i+1}, \ldots, a_{i +\ell}\} \subset A$ and $\{b_j, b_{j+1}, \ldots, b_{j +\ell}\} \subset B$ will be called an \textit{$\epsilon$-alignment} for $A$ and $B$ if $\vert a_{i+k} - b_{j+k}\vert < \epsilon$ for $k = 0,\ldots, \ell$.
\end{definition}

To find $\epsilon$-alignments, we use a simplified version of the Smith-Waterman local sequence alignment algorithm, \cite{smith1981identification}.  For the purposes of puzzle piece matching, we will search for a maximal length $\epsilon$-alignment.

%%%%%%%%%%%%%%%%%%%%%
\begin{alg} Calculate a maximum length $\epsilon$-alignment of two periodic arrays.
 \label{alg:alignment}
  \vskip 5pt
  \noindent
 \textbf{Input:}  Periodic arrays $A = \{a_0, \ldots, a_{n-1}\}$, $B = \{b_0, \ldots, b_{m-1}\}$, and $\epsilon >0$.
  \vskip 5pt 
  \noindent
  \textbf{Output:}  A maximum length $\epsilon$-alignment $\{a_i, a_{i+1}, \ldots, a_{i +\ell}\} \subset A$ and $\{b_j, b_{j+1}, \ldots, b_{j +\ell}\} \subset B$.
   \vskip 5pt
  \begin{enumerate}
    \item Construct a \textit{scoring matrix} 
   \[S_{ij} = \begin{cases}0 & \mbox{if } \vert a_i-b_j \vert <\epsilon \\ 1 & \mbox{otherwise} \\ \end{cases}\].
    \item For each $r = 0, \ldots, \gcd(m,n)$, search along the (periodic) diagonal $S_{r,0}, S_{r+1,1} \ldots$ for the longest sequence of zeros.  This longest sequence may not be unique; see Remark \ref{rem:nonunique}.
    \item Return the indices $i,j$ and length $\ell+1$, where $S_{i,j}, S_{i+1,j+1}, \ldots S_{i+\ell,j+\ell}$ is the longest sequence of zeros from 2.

\end{enumerate}
\end{alg}
%%%%%%%%%%%%%%%%%%%%%

\begin{remark}
\label{rem:nonunique}
In our application, the matrix $S_{ij}$ will have $10^5$ or more entries (e.g. comparing two arrays with lengths around $300$ or more), making it unlikely that two maximal length $\epsilon$-alignments of the same length exist.   We implement Algorithm \ref{alg:alignment} to return the most recently found maximal length $\epsilon$-alignment and do not keep track of any other $\epsilon$-alignments of the same or smaller length.
\end{remark}

Now suppose that we have puzzle pieces $P, Q$ with integral invariant signatures $A, B$ computed with a disk of radius $r$.  These pieces and their signatures are oriented counter-clockwise.  To compare puzzle pieces, we must reverse the orientation of one piece, say $Q$.  This results in a new signature $\overline {B}$ which can be obtained from $B = \{b_0, b_1, \ldots, b_{n-1}\}$ via
\[
\overline{B} = \{ \pi r^2 - b_{n-1}, \pi r^2 - b_{n-2}, \ldots, \pi r^2 - b_{1}, \pi r^2 - b_0\}.
\]
To fit $P$ and $Q$ together, we look for a maximal $\epsilon$-alignment of $A$ and $\overline B$.

\begin{definition}
Let $\epsilon > 0$.  An $\epsilon$-fit of $P = \{ p_0, \ldots, p_{m-1} \}$ with $Q=\{q_0, \ldots, q_{n-1}\}$, is a pair of substrings $\{p_i, p_{i+1}, \ldots, p_{i +\ell}\} \subset P$ and $\{q_j, q_{j+1}, \ldots, q_{j +\ell}\} \subset Q$ corresponding to an $\epsilon$-alignment $\{a_i, a_{i+1}, \ldots, a_{i +\ell}\}, \{\pi r^2 - b_{j+\ell}, \pi, \ldots, \pi r^2 - b_{j+1}, \pi r^2 - b_j\}$ of signatures $A$ and $\overline B$.
\end{definition}

\begin{remark}
Note that, although we compare $P$ and $Q$ by finding an $\epsilon$-alignment of $A$ and $\overline{B}$, the process is symmetrical:  a given $\epsilon$-alignment of $A$ and $\overline B$ corresponds to an $\epsilon$-alignment of $B$ and $\overline A$.  If the maximal $\epsilon$-fit is not unique, it is possible that the order of comparison will matter by virtue of the order in which $\epsilon$-alignments are found in Algorithm \ref{alg:alignment}.  We have not encountered this in our application.  To simplify discussion we assume that the $\epsilon$-fit of $P$ to $Q$ is the same as $Q$ to $P$.
\end{remark}

Our strategy for finding the ``best'' fit between two puzzle pieces $P$, $Q$ is to look for a maximal length $\epsilon$-fit for a well-chosen value of $\epsilon$.  This is a qualitative decision based on two main factors:
\begin{enumerate}
    \item An $\epsilon$-fit is a strict pointwise condition on the alignment of the signatures.  A typical puzzle may have short extreme changes of shape (e.g. a corner).  Our measure of fit must be sensitive to this.  Using an average measure of closeness, or allowing skips in alignment will result in incorrect fits due to these brief changes in shape (e.g. a straight edge fitting with a corner with incident straight edges).
    
    \item  A fit should be come from an alignment of maximal length.  Ideally, for standard rectangular puzzle pieces, two pieces should have an $\epsilon$-fit that includes the entire matching sides.  There will often be shorter length fits (such as the straight sides of two edge pieces or portions of edges) that come from better alignments (e.g. smaller $\epsilon$), but are not correct for puzzle assembly.
\end{enumerate}

\subsection{The orthogonal Procrustes problem}
\label{sec:procrustes}
Given an $\epsilon$-fit, the visual placement of pieces $P$ and $Q$ is done by minimizing the least squares distance between the substrings of the fit via an orientation preserving rigid motion, i.e. a transformation from $SE(2)$, the special Euclidean group.  This problem is often called the Procrustes problem,  \cite{schonemann1966generalized, eggert1997estimating}. We briefly recall the solution here for context.

Let $\{x_0, \ldots, x_\ell \}$ and $\{y_0, \ldots, y_\ell\}$ be collections of points in $\mathbb{R}^2$.  The (special) orthogonal Procrustes problem aims to find the rotation matrix $R$ that minimizes the least squares distance
\[
\sum_{i=0}^\ell \vert\vert x_i - R\, y_i\vert\vert^2.
\]
The solution to this problem is obtained via the singular value decomposition.  Viewing $x_i, y_i$ as column vectors, form the $2 \times (\ell+1)$ matrices $ X = \begin{bmatrix} x_0 & \cdots & x_\ell \end{bmatrix}$ and $Y = \begin{bmatrix} y_0 & \cdots  & y_\ell \end{bmatrix}$, and compute the singular value decomposition $XY^\top = U \Sigma V^\top$ of $X Y^\top$.  Then $R = U^\prime V^\top$, where $U^\prime$ is obtained from $U$ by multiplying the second (last) column by $\det(UV^\top)$ to ensure that $R$ has determinant one.

To apply this to our $\epsilon$-fits, we incorporate a translation to first align centroids.  Let $\{p_i, p_{i+1}, \ldots, p_{i +\ell}\}$, $\{q_j, q_{j+1}, \ldots, q_{j +\ell}\}$ be an $\epsilon$-fit of $P$ and $Q$.  Let 
\[
\displaystyle \overline p = \frac{1}{\ell+1 }\sum_{k=0}^\ell p_{i+k} \qquad \text{and} \qquad \displaystyle \overline q = \frac{1}{\ell+1}\sum_{k=0}^\ell q_{j+k}
\]
 be the respective centroids. Let
 \[
 X = \begin{bmatrix} p_i - \overline p& p_{i+1} - \overline p&\cdots & p_{i+\ell} - \overline p\end{bmatrix} \]
and
\[ Y = \begin{bmatrix} q_{j+\ell}- \overline q & q_{j+ \ell -1} - \overline q &\cdots & q_j - \overline q \end{bmatrix}
\]
 and $R$ the rotation matrix minimizing the distance between $X$ and $RY$ just described.  The transformation $g_{PQ} = (R,\overline p - R\, \overline q)$ in $SE(2)$ given by $g_{PQ}(z) = R\, z+\overline p - R\, \overline q$ then minimizes the least squares distance
 \[
\sum_{k=0}^\ell \vert\vert p_{i+k} -g \, q_{j+\ell-k}\vert\vert^2
 \]
 over all choices of $g$ in  $SE(2)$.
 In what follows, we'll write $g_{PQ}$ for the element of $SE(2)$ obtained from applying this process to an $\epsilon$-fit of $P$ and $Q$.

\begin{example}
We illustrate how an $\epsilon$-fit can vary with $\epsilon$ by applying the special Euclidean transformations obtained via Procrustes.  Shown in Figure \ref{threshold comparison} are four $\epsilon$-fits for pieces 3 and 4 from our example puzzle. These pieces have perimeters in the range $5000$ to $6500$,  and a radius $r=50$ is used for the integral invariant signatures.  For this pair of pieces, there is a large range $151 \leq \epsilon \leq 698$, for which the fit is visually correct and close to maximal length.  This behavior is typical for a correctly matched pair of rectangular puzzle pieces. 

\begin{figure*}[h]
\centering
\includegraphics[width=.85\textwidth]{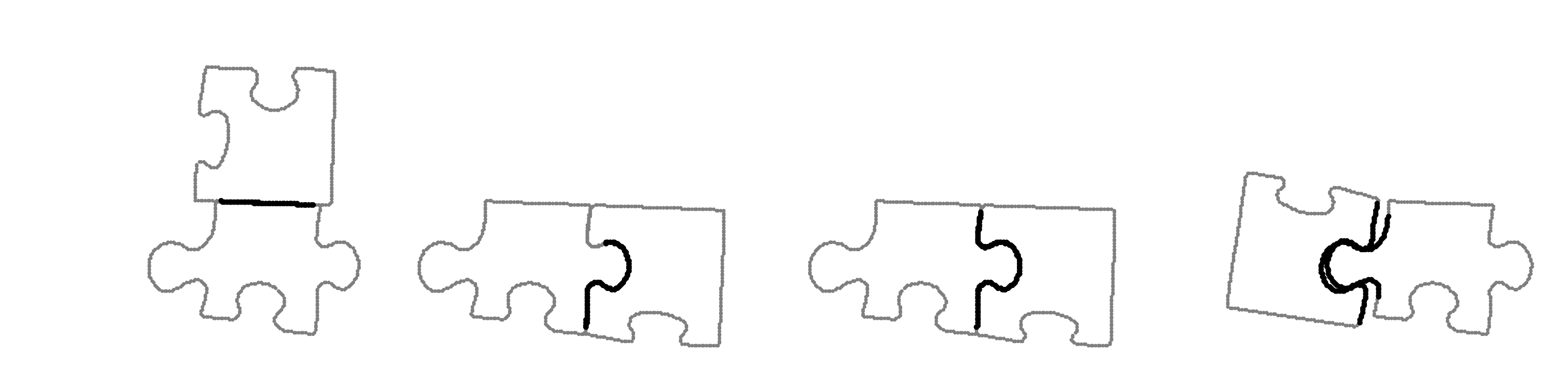}
\caption{$\epsilon$-fits with $\epsilon = 50, 150, 200,$ and $700$ (from left to right).}
\label{threshold comparison}
\end{figure*}
\end{example}

\subsection{Fit quality}

For a given $\epsilon$, there will be an $\epsilon$-fit between any pair of puzzle pieces, and it is useful to have a measure of the quality of the fit so that poor quality ones can be discarded.  Once a pairwise quality is determined, the selection of correct fits for assembly will be done globally based on this quality; this global selection will be  discussed in Section \ref{sec:assembly}.  There are two reasons we may want to discard a given $\epsilon$-fit, which we refer to as errors of type (a) and type (b):

  \begin{enumerate}
    \item[(a)] It is a fit between two pieces which \textit{do not} go together in the assembled puzzle.
    \item[(b)] It is an incorrect fit for two pieces which \textit{do} go together in the assembled puzzle.
\end{enumerate}

\begin{figure}[h]
\centering
\includegraphics[width=.35\linewidth]{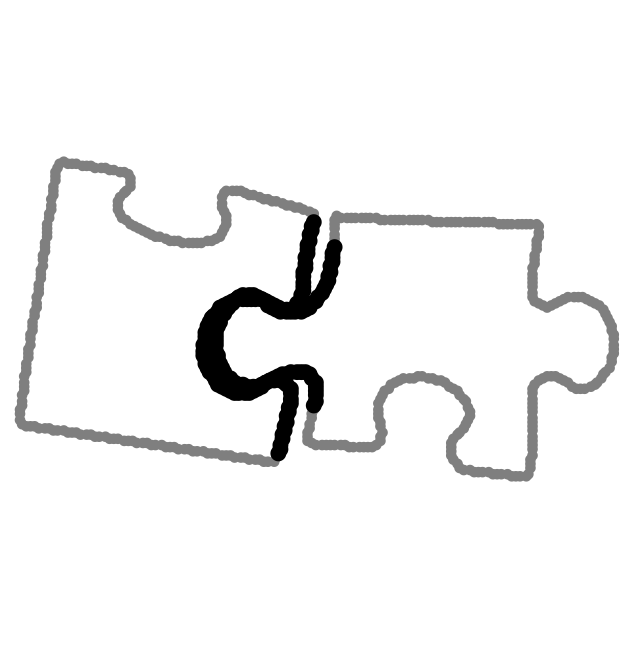}
\includegraphics[width=.35\linewidth]{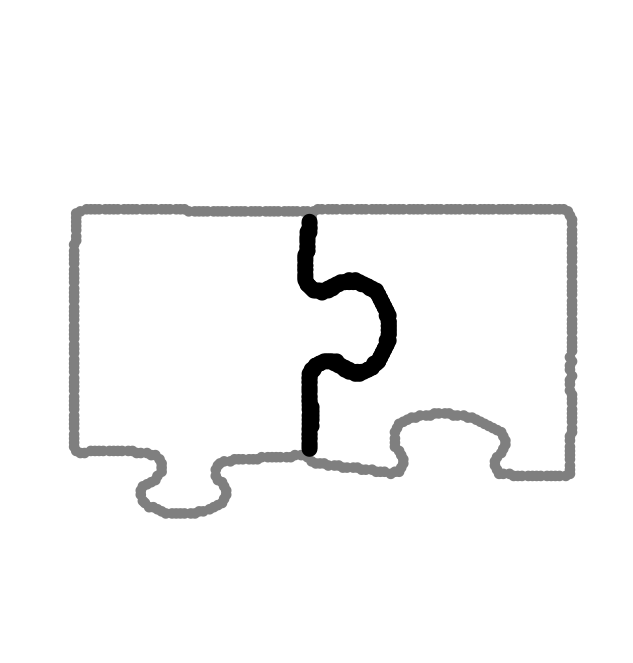}
\caption{A poor quality fit of type (b) and good quality fit of type (a).}
\label{fig:badfits}
\end{figure}

\noindent
Qualitatively, it is unlikely that one can judge whether an $\epsilon$-fit is an error of type (a) or (b) based on direct measurement of the quality of the fit.  For example, shown in Figure \ref{fig:badfits} are two $\epsilon$-fits for $\epsilon = 700$.  On the left is a clearly incorrect and poor quality fit between pieces $3$ and $4$, which do go together in the final assembly.  On the right is a good quality fit for pieces $1$ and $4$, which do not go together in the final assembly.  Despite these unavoidable flaws, some metric is needed to determine the correctness of an $\epsilon$-fit.  Our results utilize three measurements:

\begin{enumerate}
\item the length $\ell_{PQ}$ of the fit,
\item the distance
\[
d_{PQ} = \sum_{k=0}^\ell \vert\vert p_{i+k} -g_{PQ} q_{j+\ell-k}\vert\vert^2
\]
between the substrings $\{p_i, p_{i+1}, \ldots, p_{i +\ell}\}$ and $\{q_j, q_{j+1}, \ldots, q_{j +\ell}\}$ of the fit after alignment, and
\item  the standard deviations $\sigma_P, \sigma_Q$ of the substrings $\{a_i, a_{i+1}, \ldots, a_{i +\ell}\}$ and $\{b_j, b_{j+1}, \ldots, b_{j +\ell}\}$ of the invariant signatures of $P,Q$ corresponding to the fit.
\end{enumerate}

We seek a long fit with a good alignment after application of the Procrustes transformation, hence the choice of measurements 1 and 2.  Measurement 3 captures the amount of variation in the shape of the matched substrings;  small standard deviation $\sigma$ can be an indicator that our matched substrings consist of straight lines, or portions of the puzzle pieces that resemble an arc of a circle.  Our measure of the quality $q_{PQ}$ of a fit will be a function of $d_{PQ}, \ell_{PQ}$, $\sigma_{P}$ and $\sigma_{Q}$.  Other  measurements that we do not discuss here -- such as the area or perimeter of the overlap of the polygons $P$ and $g_{PQ}Q$, or the total distance between the signature substrings that give the $\epsilon$-fit -- can be effectively incorporated into the quality as well.   We use the optimization convention of minimization, so smaller will be better in our definitions of quality measurements $q_{PQ}$ of the $\epsilon$-fit of $P$ and $Q$.

The next task for assembly is to combine the $\epsilon$-fits into an assembled puzzle.  For successful assembly, correct $\epsilon$-fits need to be identified, and incorrect ones discarded.  This task is both local and global:  a fit itself can be judged based the quality measurements just discussed, while all possible fits can be considered in aggregate.  In the next section, we outline possible methods using spanning trees and the consistency of collections of $\epsilon$-fits to correctly assemble a puzzle.  

\section{Puzzle assembly}
\label{sec:assembly}

The data used to determine a successful puzzle assembly is the collection of all $\epsilon$-fits between pairs of puzzle pieces.  For a puzzle with $s$ pieces, this is a collection of $\frac{s(s-1)}{2}$ $\epsilon$-fits.  A convenient way to encode this data is as an undirected, weighted, complete graph $G$ on $s$ vertices; each vertex represents a puzzle piece, each edge an $\epsilon$-fit, and the weight on each edge a measure of the quality of the corresponding $\epsilon$-fit.  This graph $G$ will be called a \textit{comparison graph} for our puzzle.  We  briefly recall standard terminology from graph theory that will be used in this discussion.

\begin{definition}
A \textit{cycle} in a graph is a sequence of distinct edges connecting a sequence of vertices in which the only repeated vertices are the first and last in the sequence.  The \textit{length} of a cycle is the number of vertices it comprises.  We will denote a cycle by listing its vertices in order, with first and last vertex repeated.  A \textit{cycle graph} is a graph that consists of a single cycle, with no distinguished first/last vertex.  A \textit{tree} is a connected graph with no cycles.  A \textit{spanning tree} in a graph $G$ is a subgraph  of $G$ which is a tree and includes all vertices of $G$.
\end{definition}

\subsection{Spanning trees}
\begin{figure}[htbp] %  figure placement: here, top, bottom, or page
   \centering
   \includegraphics[width=\linewidth]{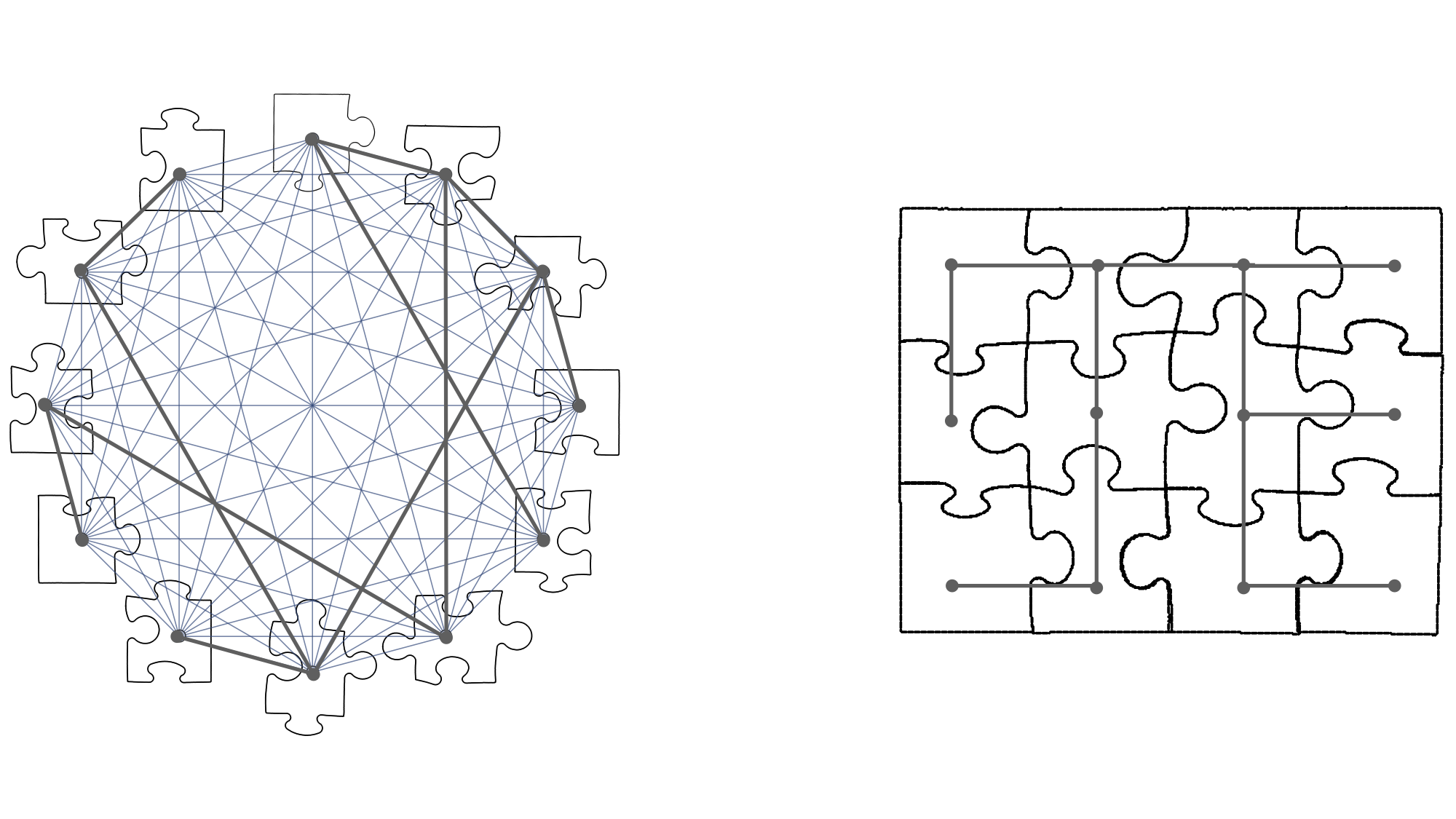} 
   \caption{A spanning tree in the complete graph $G$ and the resulting puzzle assembly.}
   \label{fig:spanningtree}
\end{figure}

In order to aggregate the collection of $\epsilon$-fits into an assembled puzzle, we choose a spanning tree in the comparison graph $G$.  This spanning tree will specify a unique way of attaching each puzzle piece to the other pieces, as shown in Figure \ref{fig:spanningtree};  if an edge connects two pieces in the spanning tree, the $\epsilon$-fit between those pieces is used in the assembly.  A spanning tree of $G$ will be called a \textit{puzzle assembly}.

To further incorporate the quality weights assigned to each edge of $G$, we choose a puzzle assembly with minimum total edge weight, which we call an \textit{optimal puzzle assembly}.  In practice, the edge weights will be unique, so the optimal puzzle assembly will also be unique.  For small puzzles for which good quality $\epsilon$-fits of type (a) and (b) are uncommon, a properly chosen measure of quality can result in optimal puzzle assemblies that are often correct or close to correct.

We illustrate various optimal puzzle assemblies using our running example puzzle.  These examples will motivate choices of quality measurement $q_{PQ}$ and the development of cycle consistency in Section \ref{sec:cycle}.  In the next two examples, we take an arclength separation $\delta = 15$ and integral invariant radius $r = 40$.   

\begin{example}
Shown in Figure \ref{fig:opa1} is the optimal puzzle assembly for quality $q_{PQ} = d_{PQ}/\ell_{PQ}$ and $\epsilon = 180$.  For this choice of $\epsilon$ there are no errors of type (b); all pieces which are meant to fit together have visually correct $\epsilon$-fits.  However, there are many good quality $\epsilon$-fits of type (a), resulting in an incorrect optimal puzzle assembly.  Many pieces are incorrectly aligned along straight edges because these offer the best quality $\epsilon$-fits for small $\epsilon$.

\begin{figure}[ht] %  figure placement: here, top, bottom, or page
   \centering
   \includegraphics[width=.4\linewidth]{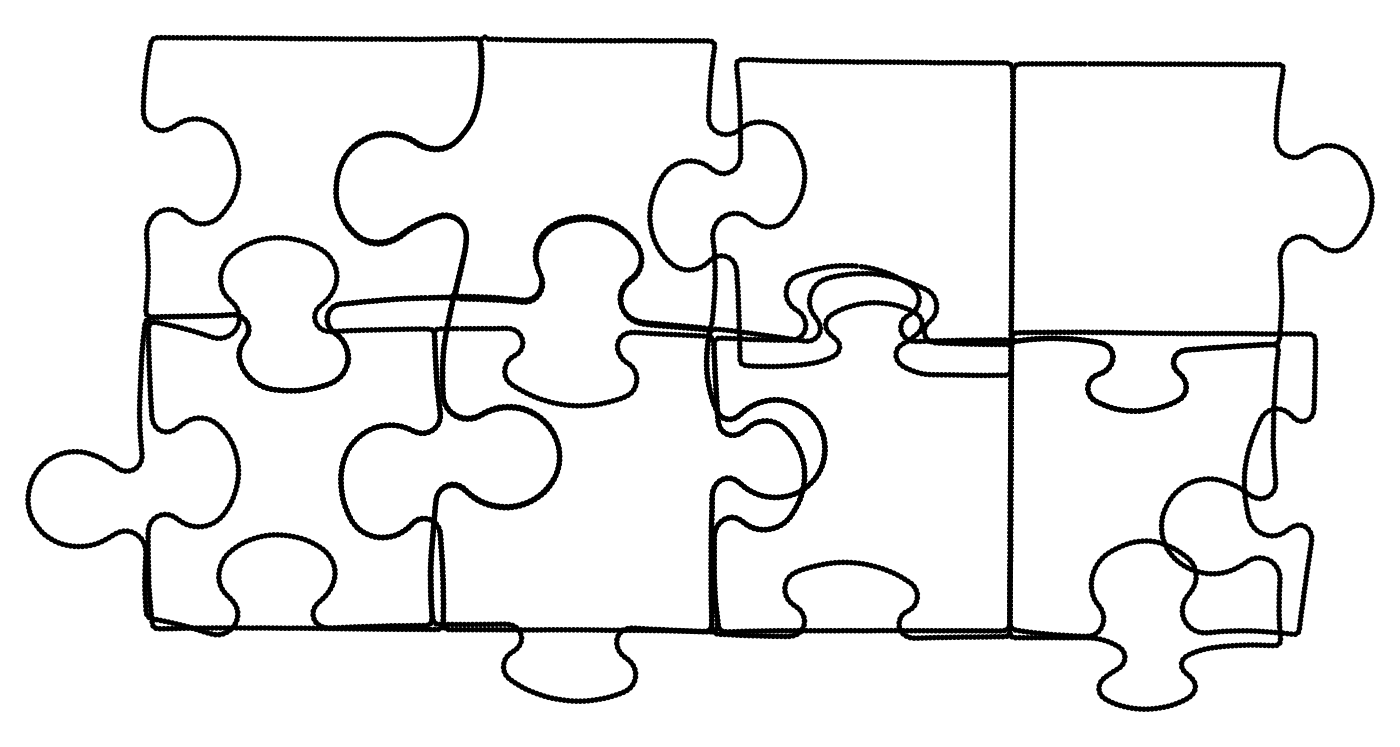}
      \includegraphics[width=.5\linewidth]{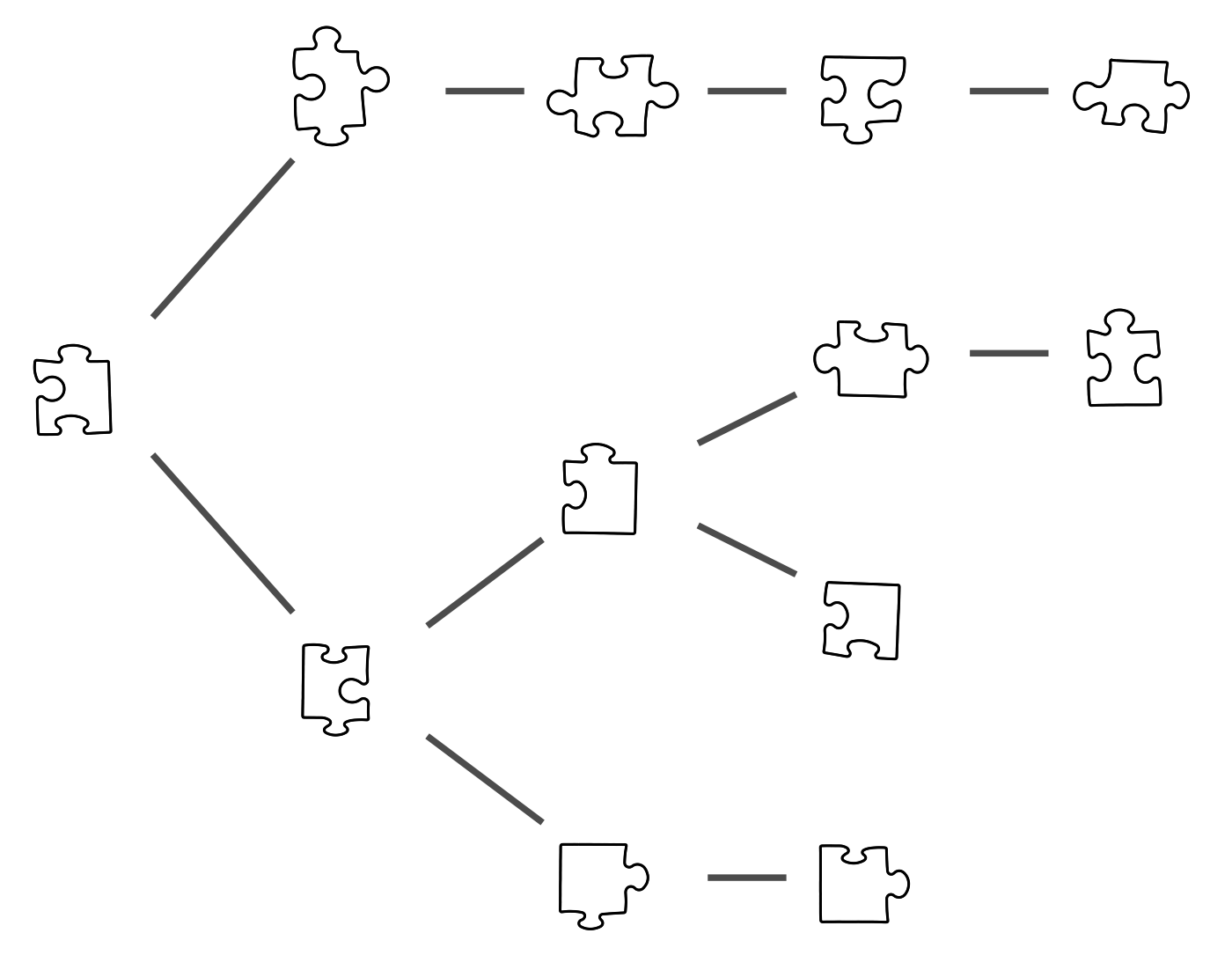} 
   \caption{An optimal puzzle assembly resulting from type (a) errors, $\epsilon = 180$.}
      \label{fig:opa1}
\end{figure}

\end{example}

It is not apparent that $d_{PQ}$ or $\ell_{PQ}$ alone can eliminate these errors of type (a), so we  incorporate $\sigma_P$ and $\sigma_Q$ via a threshold:
\begin{equation}
\label{eq:quality}
q_{PQ} = d_{PQ}/\ell_{PQ} + \iota_{\sigma^*} (\min(\sigma_P,\sigma_Q)), 
\end{equation}
where
\[
\iota_{\sigma^*} (\sigma) = \begin{cases} \infty & \mbox{if } \sigma<\sigma^* \\ 0 & \mbox{otherwise} \\ \end{cases}.
\]

A value of $\min(\sigma_P,\sigma_Q)$ that is below a chosen threshold $\sigma^*$ indicates that either substring of $P$ or $Q$ in the $\epsilon$-fit could be close to a straight line.  If this is the case, the weight $q_{PQ}$ becomes $\infty$,  eliminating the possibility that the optimal puzzle assembly will contain that particular $\epsilon$-fit. The threshold $\sigma^*$ can be chosen via a standard thresholding method.

\begin{example}
Shown in Figure \ref{fig:opa2} is the optimal puzzle assembly for quality \eqref{eq:quality} with $\epsilon = 180$ and $\sigma^* = 70$.  As before there are no errors of type (b) for this choice of $\epsilon$.  The introduction of $\sigma_P$, $\sigma_Q$ into the quality measurement now eliminates the type (a) errors resulting from straight line matches, and the resulting optimal puzzle assembly is correct.

\begin{figure}[ht] %  figure placement: here, top, bottom, or page
   \centering
   \includegraphics[width=.35\linewidth]{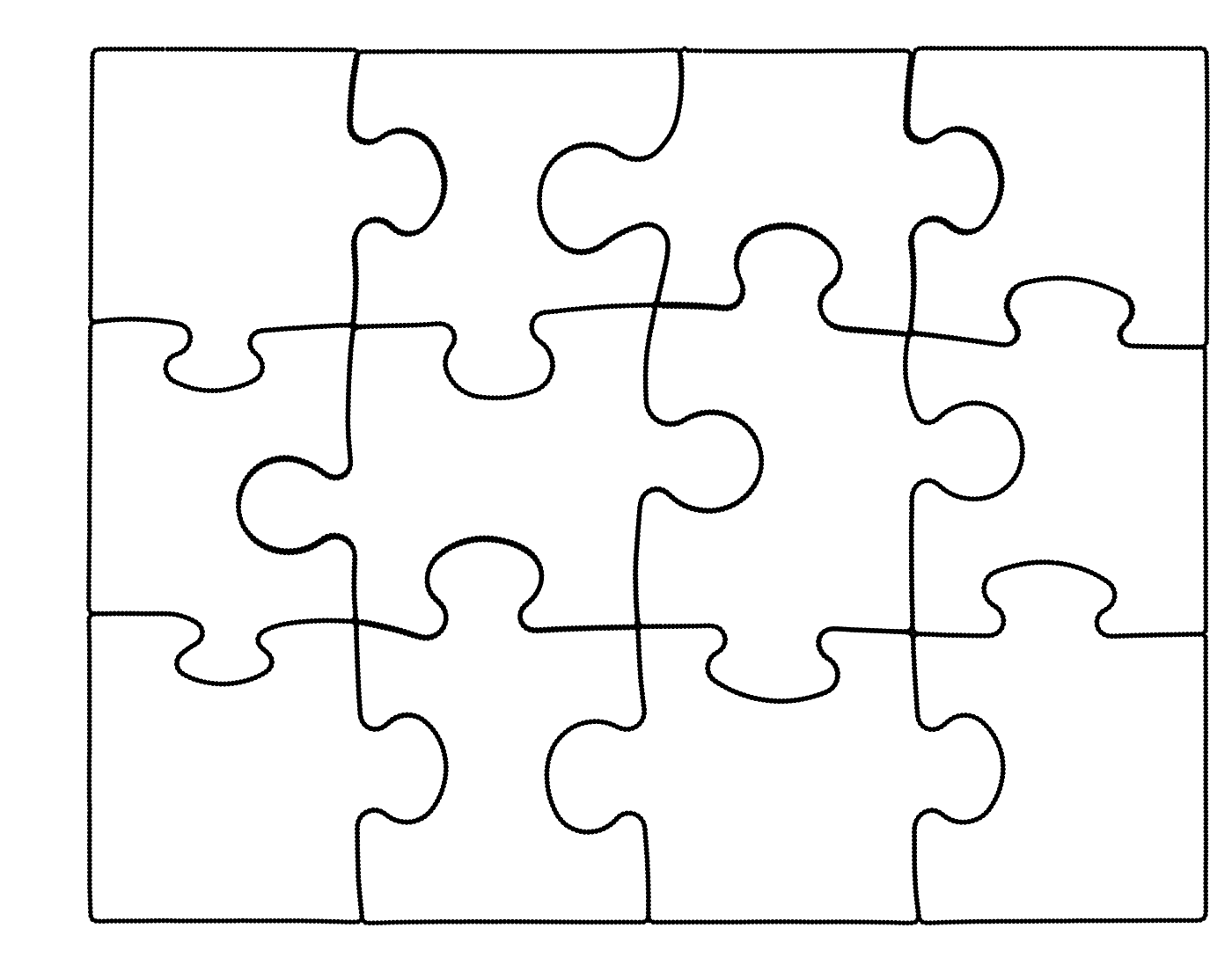}   
      \includegraphics[width=.5\linewidth]{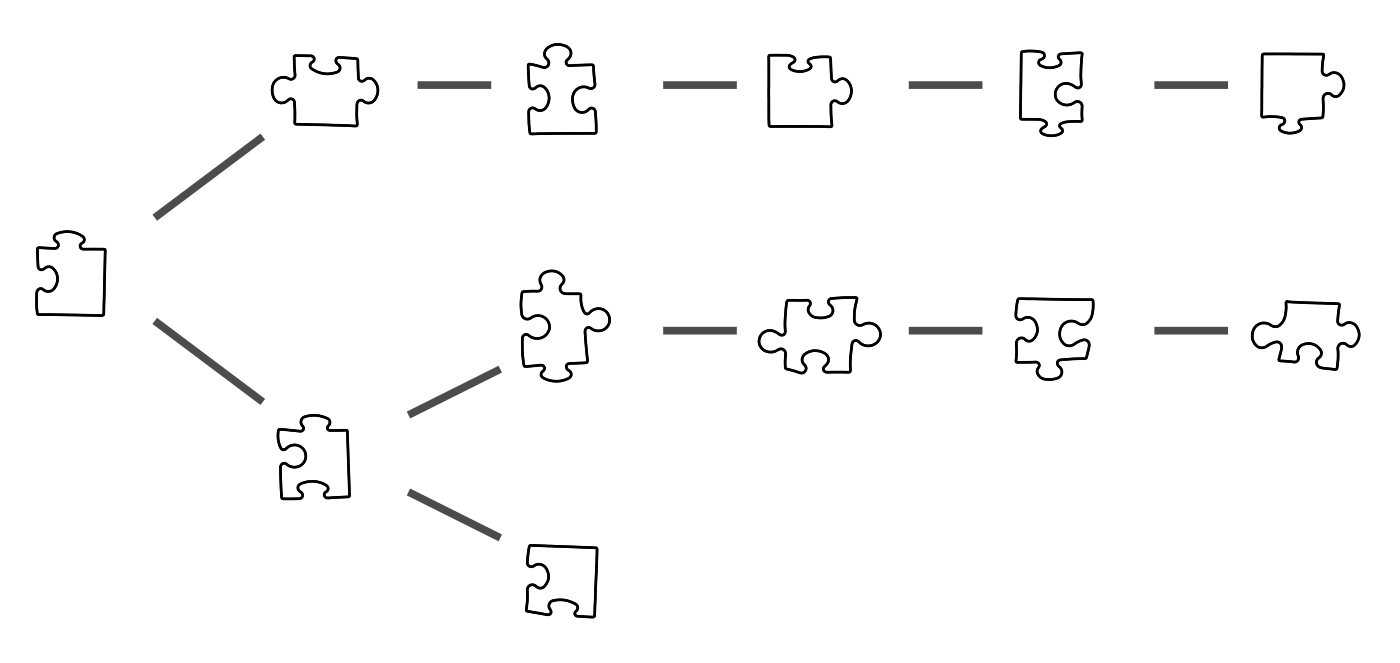} 
\caption{An optimal puzzle assembly with type (a) errors eliminated, $\epsilon = 180$.}
      \label{fig:opa2}
\end{figure}

\end{example}

\begin{example}
\label{ex:needcycleconsistency1}
Lastly, we consider another example with $\delta = 15$,  $r = 50$ and $\epsilon = 350$.  Increasing $\epsilon$ and $r$ increases the ``sloppiness'' of the fits;  there will be more potential for errors of type (a) that cannot be eliminated via the quality measurement \eqref{eq:quality}.  Shown in Figure \ref{fig:opa3}  is the optimal puzzle assembly for quality \eqref{eq:quality} with $\epsilon = 350$ and $\sigma^* = 177$.   An incorrect optimal puzzle assembly results from the type (a) errors.  In order to use this collection of $\epsilon$-fits to assemble the puzzle, we move beyond the quality measurement and examine the consistency of collections of $\epsilon$-fits; a process we will call checking \textit{cycle consistency}.

\begin{figure}[ht] %  figure placement: here, top, bottom, or page
   \centering
   \includegraphics[width=.4\linewidth]{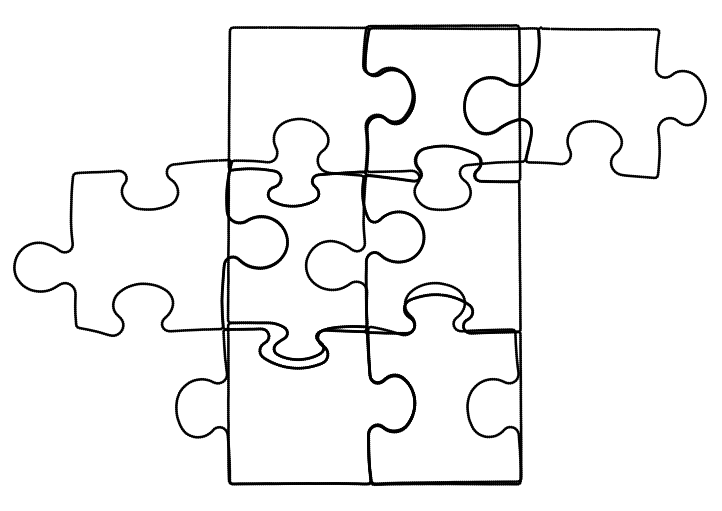}    \qquad
      \includegraphics[width=.45\linewidth]{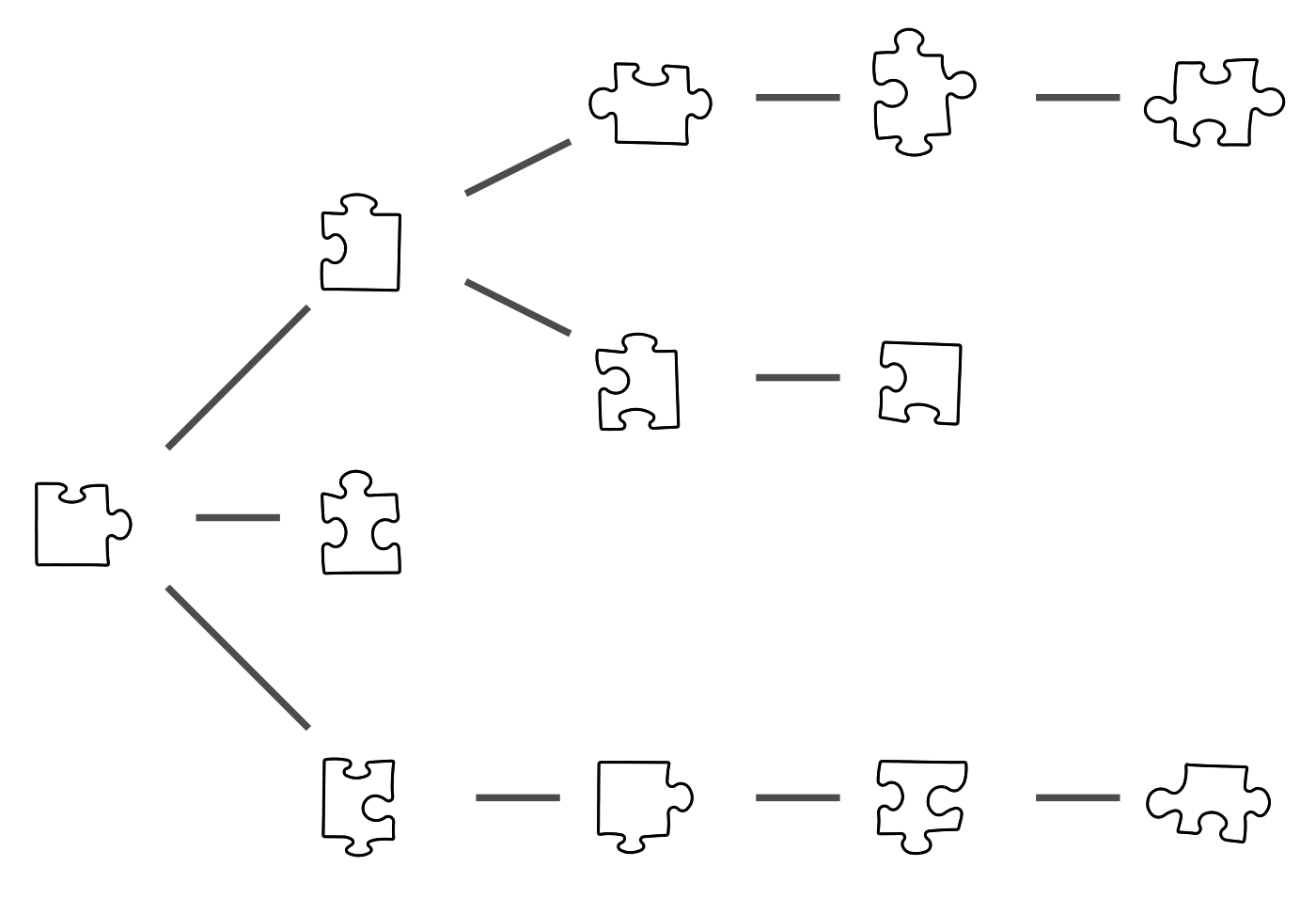} 
\caption{An incorrect optimal puzzle assembly with type (a) errors, $\epsilon = 350$.}
      \label{fig:opa3}
\end{figure}

\end{example}
%%%%%%%%%%%%%%%%%%%%%%%%
%%%%%%%%%%%%%%%%%%%%%%%%

\subsection{Cycle Consistency}
\label{sec:cycle}

Our computational solution to a jigsaw puzzle must specify the placement of each piece relative to the other pieces in a unique way, so it is mathematically convenient to interpret puzzle assemblies as spanning trees in the complete graph $G$ of all $\epsilon$-fits.  This  interpretation does not account for the fact that a puzzle has many more connections than just those chosen for the spanning tree giving the puzzle assembly.  For a standard $m \times n$ rectangular puzzle, a spanning tree uses only $mn-1$ fits, while the entire puzzle contains $2mn-(m+n)$ possible correct fits between pairs of pieces.  To use the information in these fits we propose a process of checking \textit{cycle consistency}.  

Intuitively, we would like to capture the consistency of a collection of $\epsilon$-fits.  Let $P_1, \ldots, P_s$ be the pieces of our puzzle, and write $g_{ij}$ as shorthand for the Euclidean transformation $g_{P_i P_j}$ that aligns $P_j$ with $P_i$, as discussed in Section \ref{sec:procrustes}.  Take a subcollection of $k$ pieces $P_{i_1}, \ldots, P_{i_k}$.  If, in a real puzzle, we can attach pieces in either of the sequences
\[
P_{i_1} \rightarrow P_{i_2} \rightarrow \cdots \rightarrow P_{i_{k-1}} \rightarrow P_{i_k} \quad \text{or} \quad  P_{i_1} \rightarrow P_{i_k},
\]

then the direct attachment of $P_{i_1}$ to $P_{i_k}$  should give the same placement of $P_{i_1}$ as the attachment of $P_{i_1}$ to $P_{i_k}$ through the intermediary pieces $P_{i_{2}}, \ldots, P_{i_{k-1}}$.  This collection of attachments corresponds to a cycle $(P_{i_1}, P_{i_2}, \ldots, P_{i_k},P_{i_1})$ of length $k$  the assembly graph $G$, and we propose Definitions \ref{def:transformationvalidation} and  \ref{def:overlapvalidation} as two ways to measure the consistency of the placement of pieces $P_{i_1}, \ldots, P_{i_k}$ by examining this cycle.

We first measure consistency at the level of transformations. Let $g =  g_{i_1i_2} \cdots \, g_{i_{k-1}i_k}  \, g_{i_k i_1}$
 be the composition of the  transformations in the cycle $(P_{i_1}, \ldots, P_{i_k}, P_{i_1})$.  (Note that these transformations start and end at the starting vertex $P_{i_1}$ of the cycle.) If these fits are part of a perfect puzzle assembly, $g$ should be the identity transformation.  Since our assemblies will not be perfect, we impose a threshold test on $g$.

\begin{definition}
\label{def:transformationvalidation}
Let $\theta^*, \tau^* >0$.  With $g$ as above write $g = (R_\theta, \tau)$ where $R_\theta$ is the standard rotation, $-\pi < \theta \leq \pi $, and $\tau$ the translation comprising the transformation $g$.  If $\vert \theta \vert < \theta^*$ and $\vert \vert \tau \vert \vert <\tau^*$, the cycle $(P_{i_1}, \ldots, P_{i_k},P_{i_1})$ will be called \textit{transformation consistent} (with respect to $\theta^*, \tau^*$).  If all cycles in the cycle graph $(P_{i_1}, \ldots, P_{i_k})$ are transformation consistent, then the cycle graph $(P_{i_1}, \ldots, P_{i_k})$ will be called transformation consistent.
\end{definition}

\begin{remark}
Transformation consistency for a cycle depends on the starting point of the cycle.  Thus it is possible for some cycles in a single cycle graph to be transformation consistent while others are not.  Thus checking cycle graph consistency for a cycle graph of length $k$ involves checking each of the $k$ cycles it contains.
\end{remark}

A cycle consisting of correct fits will necessarily be transformation consistent.  However, obviously incorrect fits can still be part of transformation consistent cycles, as the next example illustrates.

\begin{example}
\label{ex:tranconerr}
We return to the puzzle assembly of Example \ref{ex:needcycleconsistency1} to illustrate transformation consistency for cycles of length $4$.   There are no errors of type (b) for the chosen $\delta$, $r$ and $\epsilon$ of Example \ref{ex:needcycleconsistency1}, so we hope to discover 24 cycles corresponding to $6$ cycle graphs, shown in Figure \ref{fig:correctcycles}, that are consistent with correct assembly of the puzzle.  Removing the fits eliminated by the check on shape variation $\sigma$ (since they will not be part of the optimal puzzle assembly), we arrive at a modified comparison graph $G$ with $61$ edges and $1081$ cycle graphs of length $4$.  Using thresholds $\theta^* = \pi/20$ and $\tau^* = 30$, a total of $111$ of these cycle graphs are transformation consistent.  As shown in the example cycle of Figure \ref{fig:tranconerr}, these ``extra'' transformation consistent cycles arise from accidental alignment of $\epsilon$-fits.

\begin{figure}[ht] %  figure placement: here, top, bottom, or page
   \centering
   \includegraphics[width=.3\linewidth]{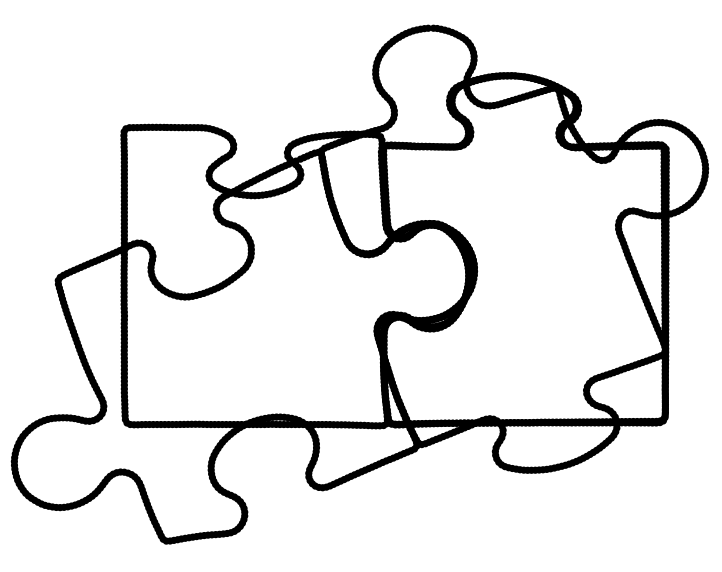}\\
   \includegraphics[width=\linewidth]{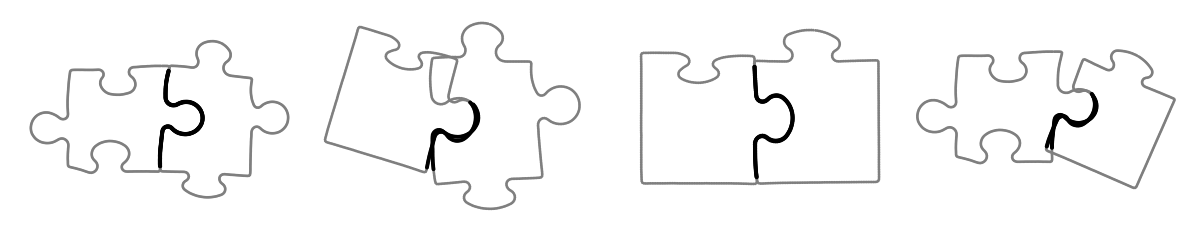} 
\caption{A transformation consistent cycle and the $\epsilon$-fits it comprises.}
      \label{fig:tranconerr} 
\end{figure}

\end{example}

Thus, transformation consistency does not provide a complete filter to identify cycles in the assembly graph consisting of correct $\epsilon$-fits.  In order to eliminate incorrect transformation consistent cycles like the one in Figure  \ref{fig:tranconerr}, we introduce \textit{overlap consistency}, a more stringent check on cycles in the comparison graph.  Overlap consistency is motivated by the simple observation that a correctly assembled cycle will not have overlapping pieces.

Denote by $\text{int}(P_i)$ the interior of the puzzle piece $P_i$ and let $P_{i_1}, \ldots, P_{i_k}$ be a collection of pieces as before.  In a perfect puzzle assembly the interiors of these pieces should not overlap as they are assembled.  That is, all pairwise intersections of the open regions
\[
\Omega_j  = \text{int}(g_{i_1 i_2}  g_{i_2 i_3} \cdots g_{i_{j-1}i_j} P_{i_j}), \quad j = 1, \ldots, k,
\]
should be empty.  Since our assemblies will not be perfect, we impose a threshold test on the areas of these intersections.  In order to account for puzzle pieces of varying size, our threshold is determined as a portion of the combined areas of the pieces being compared.

\begin{definition}
\label{def:overlapvalidation}
Let $\Omega_1, \ldots, \Omega_k$ be as above and let $\alpha^* >0$ and $\alpha^*_{ij} = \alpha^* \big( \text{area}(\Omega_i) + \text{area}(\Omega_j)   \big)$.   If $\text{area} \big( \Omega_i \cap \Omega_j \big) < \alpha^*_{ij}$ for $i,j = 1, 2, \ldots, k$, then the \textit{cycle} $(P_{i_1}, \ldots, P_{i_k},P_{i_1})$ will be called \textit{overlap consistent} (with respect to $\alpha^*$).   If all cycles in the cycle graph $P_{i_1}, \ldots, P_{i_k}$ are overlap consistent, then the \textit{cycle graph} $(P_{i_1}, \ldots, P_{i_k})$ will be called \textit{overlap consistent}.
\end{definition}

\begin{remark}
Checking overlap consistency is much more computationally intensive than transformation consistency, requiring at most $k(k-1)/2$ computations of the intersections of polygons for each  cycle.  To check overlap consistency for the whole cycle graph, these checks must be repeated for each of the $k$ starting points of cycles, naively resulting in a total of $k^2(k-1)/2$ overlap checks.  In practice, most of these overlap consistency checks will fail before checking all $k(k-1)/2$ intersections.  As needed, the computation time for overlap consistency checks can be reduced by first filtering cycles via a check on transformation consistency and by reducing the number of points in each puzzle piece by adjusting $\delta$.
\end{remark}

\begin{example}
\label{ex:overlapconsistency}
We return again to Example \ref{ex:needcycleconsistency1} and illustrate overlap consistency for cycles of length $4$.  As in Example \ref{ex:tranconerr} we use a modified comparison graph $G$ with $61$ edges and $1081$ cycles of length $4$.  With a threshold of $\alpha^* = \frac{1}{80}$, exactly $6$ of these cycle graphs are overlap consistent, shown in Figure  \ref{fig:correctcycles}.  These are the $6$ correct cycle graphs of length $4$ that we can expect to find in a perfect assembly. 

\begin{figure*}[ht] %  figure placement: here, top, bottom, or page
   \centering
 \includegraphics[width=.9\linewidth]{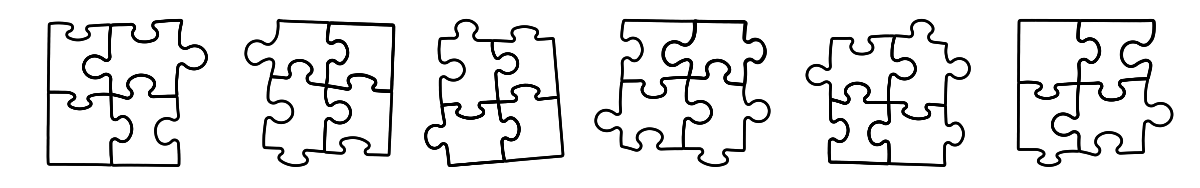}
\caption{Overlap consistent cycles for the $\epsilon$-fits of Example \ref{ex:needcycleconsistency1} with $\alpha^* = \frac{1}{80}$.}
      \label{fig:correctcycles}
\end{figure*}
\end{example}

To improve our optimal puzzle assembly, we incorporate cycle consistency information into the edge weights of the assembly graph.  We do this by reducing the edge weights for those fits which appear in overlap consistent cycle graphs, increasing the likelihood that these fits are selected in the optimal puzzle assembly.  

\begin{definition}
Suppose that the $\epsilon$-fit of $P_i$ and $P_j$ has weight $q_{ij}$ and appears in $c$ overlap consistent cycles.  Choose $0 < \beta^* < 1$, and assign a new weight $ \overline{q_{ij}} = (\beta^*)^c q_{ij}$ to the $\epsilon$-fit.  For rectangular puzzles, we will have $c=0,1,$ or $2$.  The graph $\overline{G}$ obtained from $G$ by assigning the new weights $\overline{q_{ij}}$ will be called the \textit{cycle consistent comparison graph}.
\end{definition}

\begin{example}
\label{ex:needcycleconsistency3}
We return to Example \ref{ex:needcycleconsistency1} for a final time and incorporate cycle consistency into the optimal puzzle assembly.  Using the overlap consistent cycle graphs found in Example \ref{ex:overlapconsistency} and taking $\beta^* = \frac{1}{2}$, we create the new overlap consistent comparison graph $\overline G$.  Taking a minimal spanning tree of this graph now results in a correct optimal puzzle assembly, shown in Figure \ref{fig:opa4}.  Hence cycle consistency is able to discern the correct $\epsilon$-fits to include in the assembly where the quality measure alone is not.

\begin{figure}[ht] %  figure placement: here, top, bottom, or page
   \centering
 \includegraphics[width=.35\linewidth]{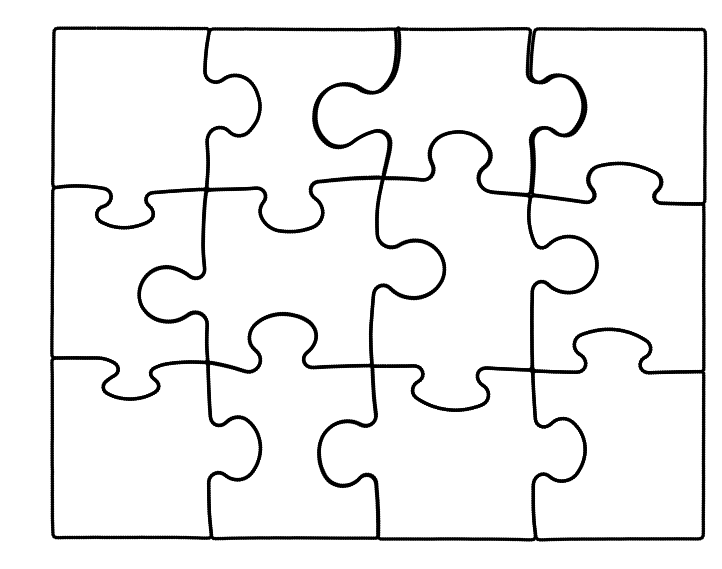}    \qquad
      \includegraphics[width=.45\linewidth]{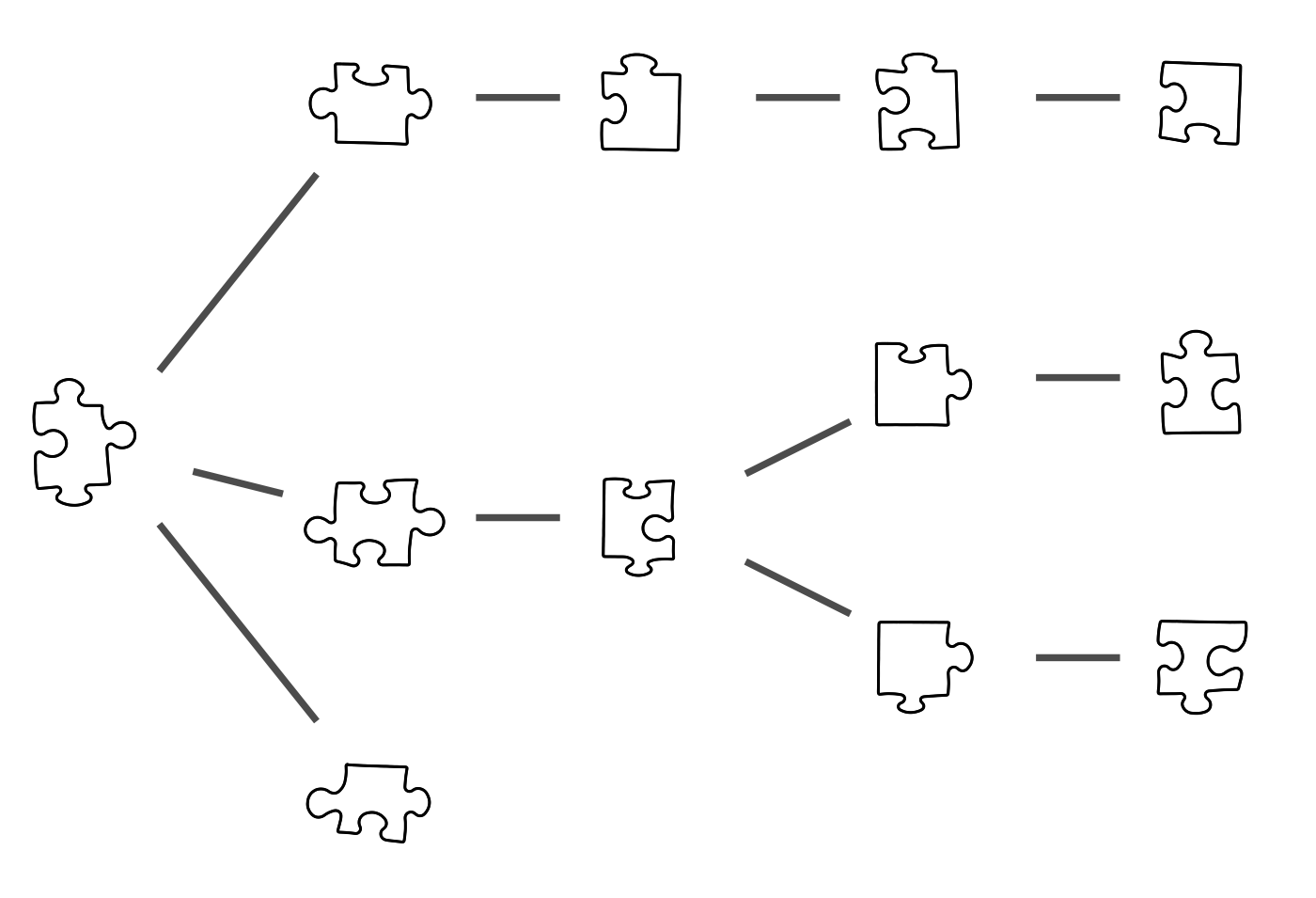} 
\caption{An optimal puzzle assembly incorporating cycle consistency, $\epsilon = 350$.}
      \label{fig:opa4}
\end{figure}
\end{example}

\subsection{The full puzzle assembly algorithm}

We summarize our full method for automated puzzle assembly.  Assembly results from applying this method to various example puzzles will be provided in Section \ref{sec:results}.

%%%%%%%%%%%%%%%%%%%%%
\begin{alg}  An algorithm for apictorial jigsaw puzzle assembly.
\label{alg:fullassembly}
  \vskip 5pt
  \noindent
 \textbf{Input:}  A collection of (unprocessed) puzzle piece boundary curves.
  \vskip 5pt 
  \noindent
  \textbf{Output:}  A collection of transformations intended to produce a correct puzzle assembly.
   \vskip 5pt
  \begin{enumerate}
    \item \textit{Process the puzzle data.}  Choose an arclength separation $\delta$ and apply Algorithm \ref{alg:resampling} to the unprocessed pieces until arclength separation is uniform.  In practice, this takes fewer than $10$ iterations.  Let $P_1, \ldots, P_s$ be the processed pieces.
    
    \item \textit{Compute integral invariants.}  Choose an integral invariant radius $r$ and apply Algorithm  \ref{alg:areainvariant} to compute discrete integral area invariants $A_1, \ldots, A_s$ for $P_1, \ldots, P_s$.
    
    \item \textit{Compute $\epsilon$-fits for all piece pairs.}  Choose a fit threshold $\epsilon$ and apply Algorithm  \ref{alg:alignment} to find the maximum length $\epsilon$-fit $P_i \leftarrow P_j$ for each pair of pieces $P_i,P_j$.  Store the quality data $\ell_{ij}$, $d_{ij}$, $\sigma_{i}$, $\sigma_{j}$ and the Procrustres transformation $g_{ij}$ mapping $P_j$ to $P_i$.
    
     \item \textit{Form the comparison graph.} Choose a fit quality $q_{ij} = f(\ell_{ij}, d_{ij}, \sigma_i, \sigma_j)$  and form the weighted comparison graph $G$ as a complete graph edges $P_i$ and edge weights $q_{ij}$.  Optionally, the minimal spanning tree of $G$ can be computed and tested to see if it yields a correct optimal puzzle assembly.
     
     \item[5a.] \textit{Check cycle transformation consistency (optional, as a prefilter to overlap consistency).}  Choose thresholds $\theta^*, \tau^*$ and  a collection of cycles in $G$ to check for transformation consistency with respect to $\theta^*, \tau^*$.
     
      \item[5b.] \textit{Check cycle overlap  consistency (optional, if comparison graph does not produce correct assembly).}  Choose a threshold $\alpha^*$ and a collection of cycles in $G$ (e.g. the transformation consistent cycles from 5) to check for overlap consistency with respect to $\alpha^*$.
     
     \item[6.] \textit{Form the cycle consistent comparison graph.} With $0 <  \beta^* < 1$, adjust the weights of $G$ to create the cycle consistent comparison graph $\overline G$.  The minimal spanning tree of $\overline G$ can then be tested to see if it yields a correct optimal puzzle assembly.
\end{enumerate}
\end{alg}

%%%%%%%%%%%%%%%%%%%%%

\section{Results}
\label{sec:results}

\begin{example}
\label{ex:alphabet}
Pictured in Figure \ref{fig:alphabet} is the assembly of the 50 piece ``Alphabet'' puzzle,  \cite{alphabet}.  Applying Algorithm \ref{alg:fullassembly} with $\delta = 15, r= 50, \epsilon = 220, \sigma^* = 115$ and quality 
\[
q_{ij} = d_{ij}/\ell_{ij} + \iota_{\sigma^*} (\min(\sigma_i,\sigma_j))
\]
 yields this puzzle assembly.  The optimal puzzle assembly obtained directly from $G$ is sufficient to solve this puzzle; no cycle consistency checks were performed.  The entire computation (including segmentation of the puzzle images) required about 3 minutes on an Apple Macbook Air M1 8gb using Mathematica 12.3.1.0.

\begin{figure*}[ht] %  figure placement: here, top, bottom, or page
   \centering
   \includegraphics[width=.5\linewidth]{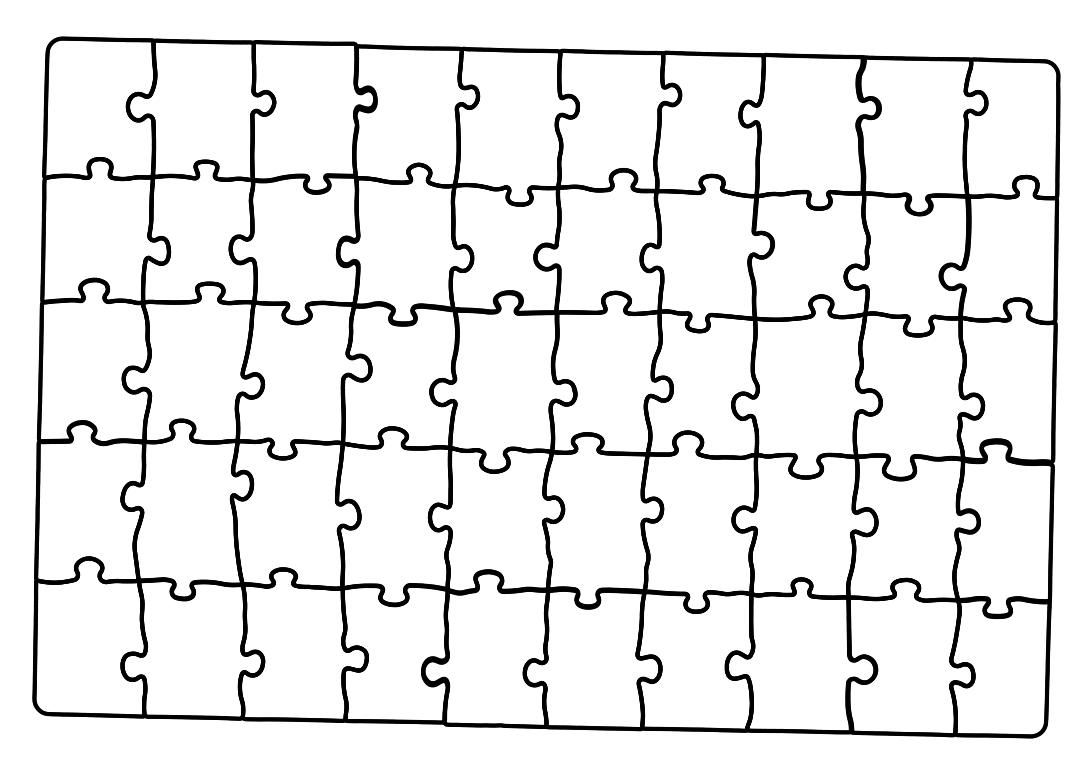} 
\caption{An optimal puzzle assembly for the Alphabet puzzle}
      \label{fig:alphabet}
\end{figure*}
\end{example}

\begin{example}
Pictured in Figure \ref{fig:rainforest} is the assembly of the 46 piece ``Rainforest'' puzzle,  \cite{rainforest}.  Applying Algorithm \ref{alg:fullassembly} with the same $\delta, r, \epsilon,  \sigma^*$ and $q_{ij}$ as in Example  \ref{ex:alphabet} yields this puzzle assembly.  The optimal puzzle assembly obtained directly from $G$ is sufficient to solve this puzzle; no cycle consistency checks were performed.  The entire computation (including segmentation of the puzzle images) required about 3 minutes on an Apple Macbook Air M1 8gb using \texttt{Mathematica} 12.3.1.0.

\begin{figure*}[ht] %  figure placement: here, top, bottom, or page
   \centering
 \includegraphics[width=.7\linewidth]{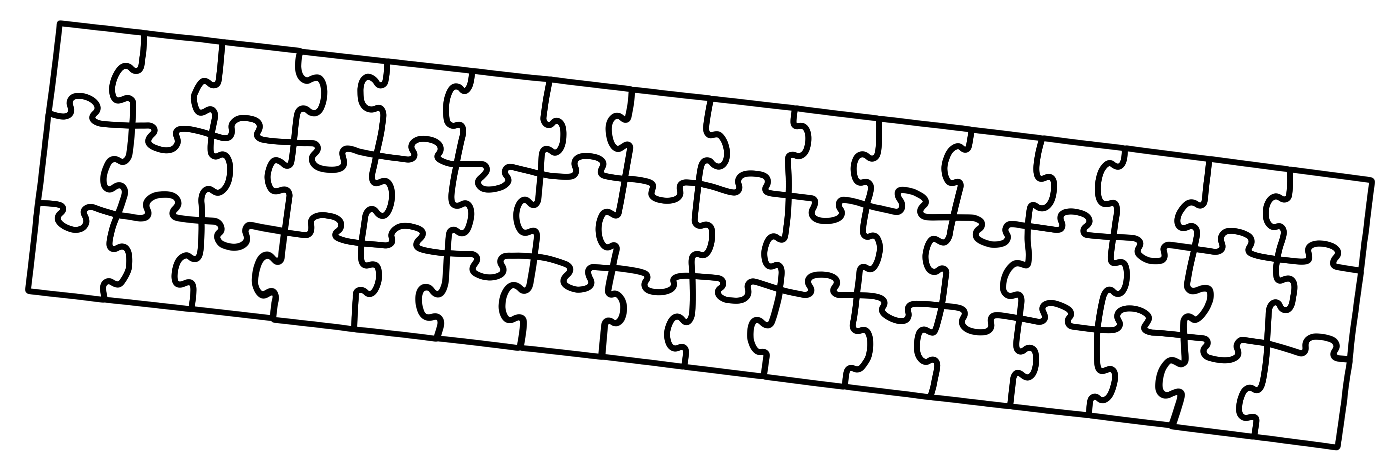}    \qquad
\caption{An optimal puzzle assembly for the Rainforest puzzle.}
      \label{fig:rainforest}
\end{figure*}

\end{example}

\begin{example}
\label{ex:safari}
Pictured in Figure \ref{fig:safari1} are two assemblies of the 100 piece ``Safari'' puzzle,  \cite{safari}.  For both assemblies, we use $\delta = 15, r= 50, \epsilon = 260, \sigma^* = 205$ and the same quality $q_{ij}$ as in the previous examples.  For the first assembly, no cycle consistency checks were performed, and the assembly is very incorrect.  This assembly required  9 minutes Apple Macbook Air M1 8gb using \texttt{Mathematica} 12.3.1.0.  For the second assembly,  two stages of cycle consistency checks with cycle length 4 are performed:  first transformation consistency with $\theta^* = \pi/30$ and $d^* = 20$, then overlap consistency with $\alpha^* = 1/80$ on the transformation consistent cycles.  In the comparison graph there are $9057744$ cycles of length 4, so this computation is intensive, requiring roughly 3 additional hours for cycle consistency checks.  Note that the assembly is much improved, with only one piece misplaced.

\begin{figure*}[ht] %  figure placement: here, top, bottom, or page
   \centering
 \includegraphics[width=.6\linewidth]{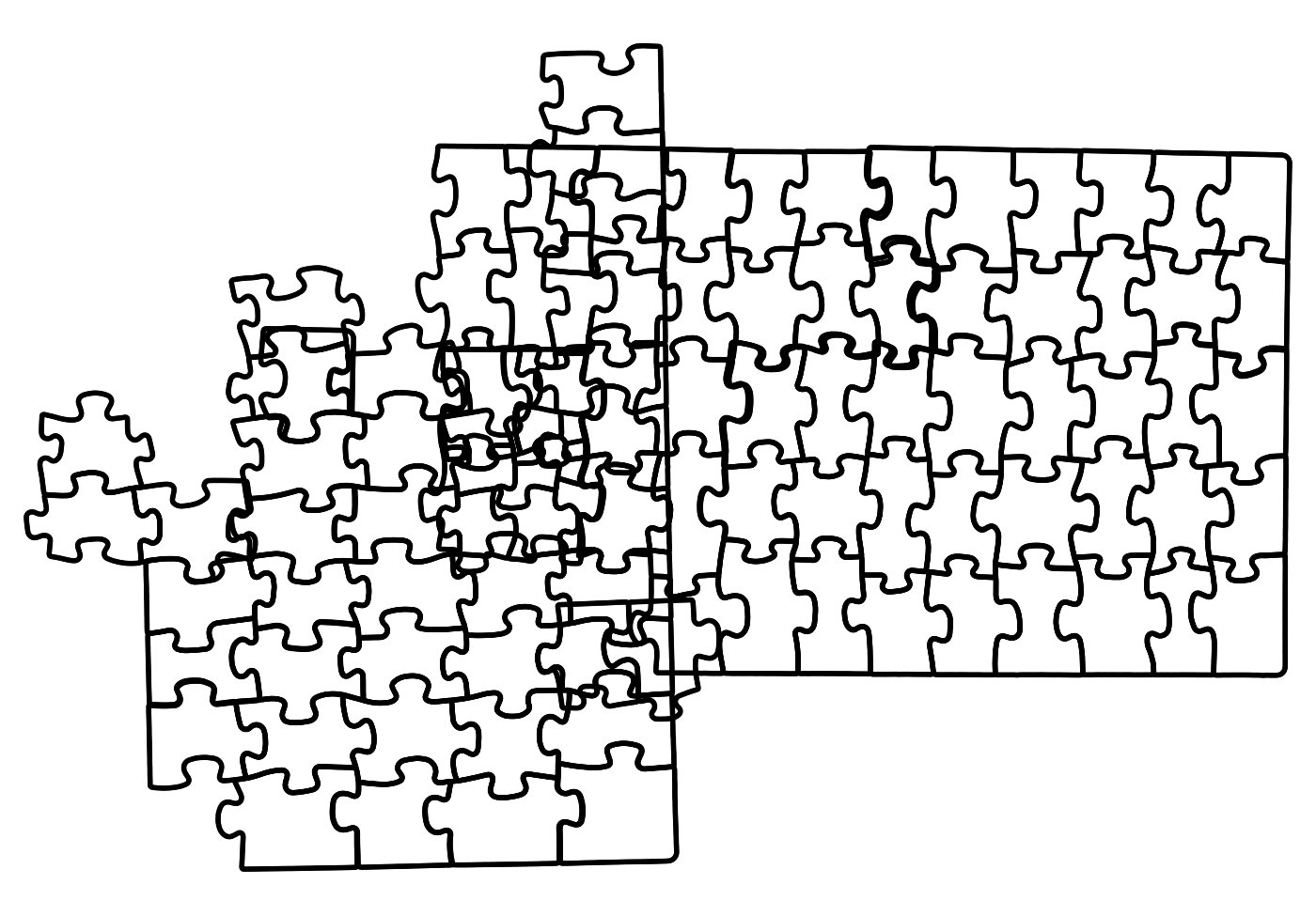}    \qquad
      \includegraphics[width=.75\linewidth]{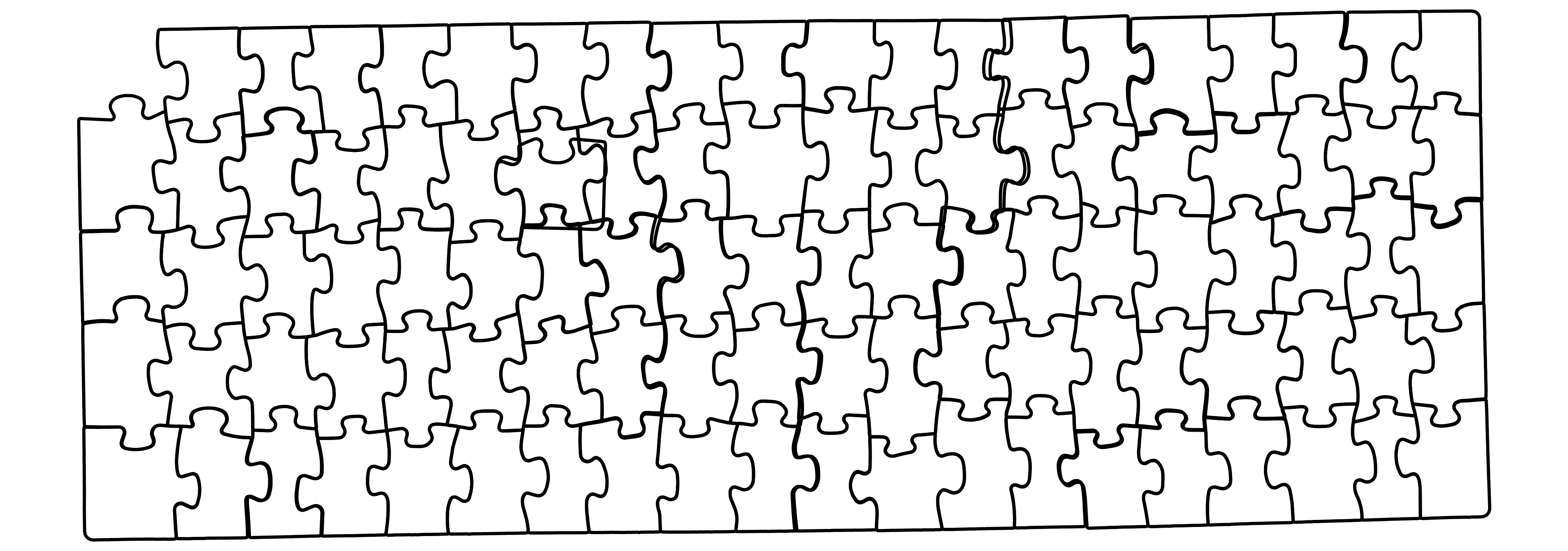} 
\caption{Optimal assemblies for the Safari puzzle, with and without cycle consistency.}
      \label{fig:safari1}
\end{figure*}

\begin{figure*}[ht] %  figure placement: here, top, bottom, or page
   \centering
      \includegraphics[width=.75\linewidth]{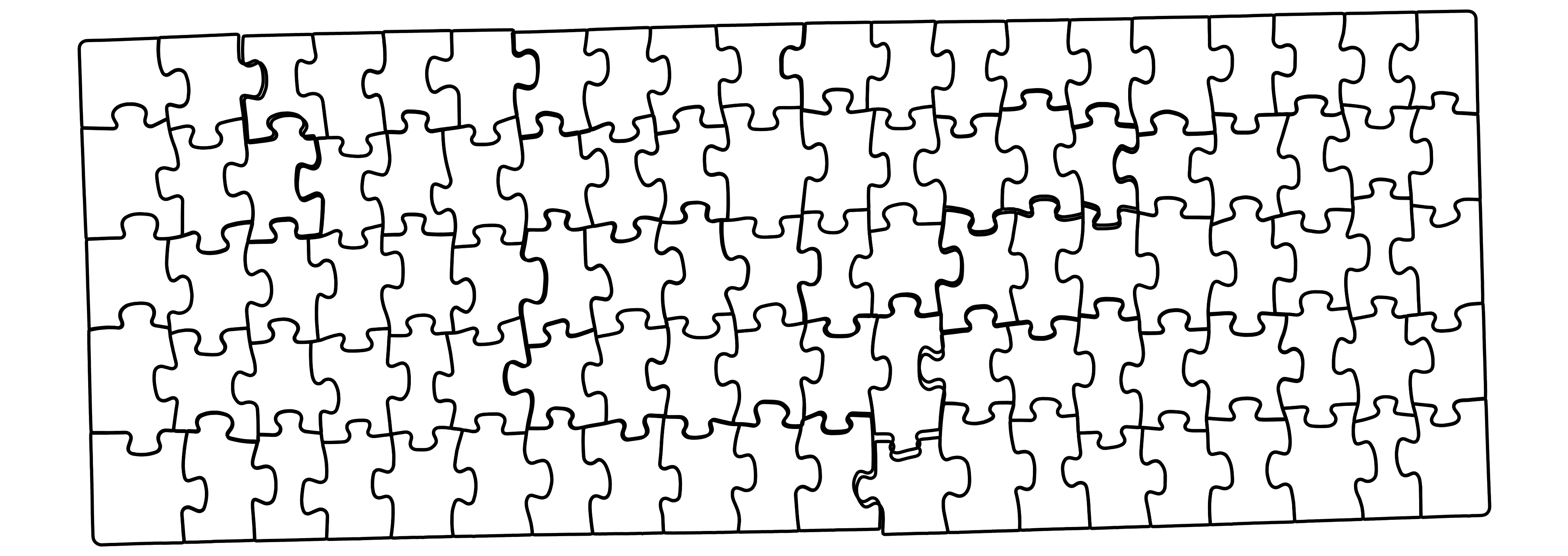} 
\caption{Optimal assembly for the Safari puzzle, with a modified quality $q_{ij}$.}
      \label{fig:safari2}
\end{figure*}

Surprisingly, with the same values of $\delta, r, \epsilon$ and  $\sigma^*$, we can achieve a completely correct puzzle assembly by changing the quality to 
\[
q_{ij} = d_{ij}/\ell_{ij}^3 + \iota_{\sigma^*} (\min(\sigma_i,\sigma_j)).
\]
This change puts stronger weight on the length of the match, rather than just the average distance between paired points in the $\epsilon$-fit.  Shown in Figure \ref{fig:safari2} is the (correct) assembly for this quality, with no cycle consistency checks performed.  This assembly required  9 minutes on the same hardware.
\end{example}

\section{Conclusion}

We have described a new method for automatic assembly of jigsaw puzzles.  This method highlights the efficacy of integral area invariants for shape comparison, and uses combinatorial information about the collection of shape matches to better select which matches should be used for the puzzle assembly.  As shown in Section \ref{sec:results}, our method is effective at assembling traditional rectangular jigsaw puzzles, but does not rely on structural information about piece shape or arrangement.  There are a number of interesting directions for further research.

In \cite{hoff2014automatic}, related methods for puzzle assembly are shown to be effective on non-rectangular puzzles such as the Baffler Nonagon, \cite{baffler}.  This robustness to irregular assemblies is helpful for solving more general object reassembly problems.  For reasons that are not yet clear, preliminary experiments applying our assembly method to non-rectangular puzzles have not yielded results as impressive as those of \cite{hoff2014automatic}.  It would be worthwhile to better understand this shortcoming, and to experiment with changes to Algorithm \ref{alg:fullassembly} to make it more effective on non-rectangular puzzles.

Our method of determining $\epsilon$-fits involves exhaustive comparison of  invariant signatures.  This was a qualitative choice based on the goal of finding the longest match within a threshold of shape similarity.  Other comparison methods or paradigms for comparison could be explored.  Rather than enforcing a strict distance threshold, signature matches could be found using a more flexible notion of local sequence alignment, e.g. one that measures distance with allowance for errors or omissions, \cite{navarro2001guided}.  One could also search for alignments probabilistically, \cite{eddy2008probabilistic}.  Additionally, our shape comparisons were done independently; the fit of one pair of pieces does not affect the fit of another pair.  This does not comport with reality, where a correct fit between a pair of pieces limits the further possible fits for those pieces.  To reflect this fact in the algorithmic approach it may be possible to use multiple sequence alignment, \cite{feng1987progressive}, which could seek non-overlapping alignments of a collection of signatures.

For large puzzles and large cycles, the cycle consistency process becomes computationally intensive since the number of cycles grows exponentially with both the number of vertices and the length of the cycle.  Optimization of the approach to cycle consistency is needed to use it effectively to assemble puzzles larger than 100 pieces.  One possible approach to this optimization would be to select cycles randomly to test, and update the cycle consistent comparison graph $\overline G$ dynamically.  Random cycle selection could be done in an informed manner; prioritizing those consisting of $\epsilon$-fits with better quality, for example.

There are a number of parameters that need to be determined for a successful puzzle assembly.   Some of these parameters, such as the arclength resolution $\delta$, the integral invariant radius $r$ and the comparison threshold $\epsilon$,  seem to depend primarily on the size of the puzzle pieces and not on the puzzle geometry \textit{per se}.  Others, such as the shape variation threshold $\sigma^*$ and the overlap threshold $\alpha^*$ depend on the geometry and arrangments of the pieces.  It would be worthwhile to investigate how some of these parameters could be determined or optimized automatically.  Along similar lines, one could investigate the relationship between parameter values and quality measurements of the fits.  As can be seen in Example \ref{ex:safari}, the choice of quality can have a strong effect on the correctness of the assembly, and a more systematic method for choosing quality would be useful.

Finally,  because the focus of this work was on the application of integral area invariants, the use of other invariant signatures (e.g. differential invariants, \cite{calabi1998differential,boutin2000numerically}, or invariant histograms, \cite{brinkman2012invariant}) in combination with our graph based assembly process was not systematically explored.  It would be very interesting if other types of invariant signatures proved to be more (or less) effective in achieving correct puzzle assemblies.

%%%%%%%%%%%%%%%%%%%%%%%%
%%%%%%%%%%%%%%%%%%%%%%%%
%%%%%%%%%%%%%%%%%%%%%%%%

\section{Acknowledgements}

The Carleton College Towsley Endowment provided funding that made this research possible.  We would like to thank Peter Olver and Irina Kogan for helpful conversations.  We are also grateful to Noah Goldman, Marshall Ma, and Jason Zhu for their contributions to an early version of this project.

\section{Data Availability}
Selected \texttt{Mathematica} code and puzzle image data is available at the second author's website, 
\url{http://people.carleton.edu/~rthompson}.
%\[
%\text{{\tiny
%\url{https://drive.google.com/drive/folders/16xeJk5VTViYYWRGcY4IvM2n7AKTycdTS?usp=sharing}}}
%\]
Additional code and data  may be requested from the second author.
\newpage

\bibliographystyle{unsrt}

\bibliography{references}% common bib file


\end{document}